\documentclass[a4paper,10pt,,reqno]{amsart}
\usepackage[utf8]{inputenc}
\usepackage{amssymb}
\usepackage[dvips]{graphicx}

\usepackage{color}

\usepackage[colorlinks=true]{hyperref}

\definecolor{dark-red}{rgb}{.54,.0,.0}
\definecolor{dark-green}{rgb}{.0,.4,.0}
\definecolor{dark-blue}{rgb}{.04,.04,.4}

\hypersetup{linkcolor=dark-red, urlcolor=dark-blue, citecolor=dark-green}

\newcommand{\cal}{\mathcal}

\begin{document}

\newcommand{\lsu}{Lady\v{z}enskaja, Solonnikov, and Ural'ceva}
\newcommand{\tch}{Taylor, Cahn, and Handwerker}
\newcommand{\thom}{Thom\'{e}e}
\newcommand{\cSection}{Section}
\newcommand{\cSections}{Sections}
\newcommand{\csection}{section}
\newcommand{\csections}{sections}
\newcommand{\cpaper}{paper}
\newcounter{transfer}

\newcounter{growth}
\newcounter{ldot}
\newcounter{keep}
\newcounter{jump}
\newcounter{maximum}
\newcounter{ratio}
\newcounter{curvature}

\newcommand{\half}{\mbox{$\displaystyle\frac{1}{2}$}}
\newcommand{\onesixth}{\mbox{$\displaystyle\frac{1}{6}$}}
\newcommand{\onethird}{\mbox{$\displaystyle\frac{1}{3}$}}
\newcommand{\onefourth}{\mbox{$\displaystyle\frac{1}{4}$}}
\newcommand{\ux}{\mbox{$u_x^2$}}
\newcommand{\ui}{\mbox{$(u_x^m)_i$}}
\newcommand{\up}{\mbox{$(u_x^m)_{i+1}$}}
\newcommand{\upp}{\mbox{$(u_x^m)_{i+2}$}}
\newcommand{\um}{\mbox{$(u_x^m)_{i-1}$}}
\newcommand{\Wi}{\mbox{$W[(u_x^m)_i]$}}
\newcommand{\Wp}{\mbox{$W[(u_x^m)_{i+1}]$}}
\newcommand{\Wm}{\mbox{$W[(u_x^m)_{i-1}]$}}
\newcommand{\uti}{\mbox{$(u_x^m)_i^2$}}
\newcommand{\utp}{\mbox{$(u_x^m)_{i+1}^2$}}
\newcommand{\utm}{\mbox{$(u_x^m)_{i-1}^2$}}
\newcommand{\pWi}{\mbox{$\displaystyle \frac{\Wp\ -\ \Wi}
                        {\up\ -\ \ui}$}}
\newcommand{\pWm}{\mbox{$\displaystyle \frac{\Wi\ -\ \Wm}
                        {\ui\ -\ \um}$}}
\newcommand{\ppW}{\mbox{$\pWi\ -\ \pWm$}}
\newcommand{\WD}{\mbox{$W^{\prime\prime}$}}
\newcommand{\dW}{\mbox{$W^{\prime}$}}
\newcommand{\Et}{\mbox{$\bar{E}$}}
\newcommand{\Wt}{\mbox{$\bar{W}$}}
\newcommand{\Wb}{\mbox{$\tilde{W}$}}
\newcommand{\Wtp}{\mbox{$\Wb^{\prime}$}}
\newcommand{\bE}{\mbox{$E_i$}}
\newcommand{\Ehat}{\mbox{$\hat{E}_i$}}
\newcommand{\Ehats}{\mbox{$\hat{E}_i^2$}}
\newcommand{\EEE}{\mbox{$E$}}
\newcommand{\ci}{\mbox{$c_i^0$}}
\newcommand{\cp}{\mbox{$c_i^1$}}
\newcommand{\cm}{\mbox{$c_i^{-1}$}}
\newcommand{\tui}{\mbox{$(u_t^m)_i$}}
\newcommand{\tuis}{\mbox{$(u_t^m)_i^2$}}
\newcommand{\tup}{\mbox{$(u_t^m)_{i+1}$}}
\newcommand{\tum}{\mbox{$(u_t^m)_{i-1}$}}
\newcommand{\ddt}{\mbox{$\displaystyle\frac{d}{dt}$}}
\newcommand{\xid}{\mbox{$\dot{x}_i$}}
\newcommand{\xpd}{\mbox{$\dot{x}_{i+1}$}}
\newcommand{\xmd}{\mbox{$\dot{x}_{i-1}$}}
\newcommand{\upmum}{\mbox{$\up\ -\ \um$}}
\newcommand{\uxi}{\mbox{$u(x_i(t),t)$}}
\newcommand{\uxp}{\mbox{$u(x_{i+1}(t),t)$}}
\newcommand{\uxm}{\mbox{$u(x_{i-1}(t),t)$}}
\newcommand{\uxxi}{\mbox{$u_x(x_i(t),t)$}}
\newcommand{\uxxp}{\mbox{$u_x(x_{i+1}(t),t)$}}
\newcommand{\uxxm}{\mbox{$u_x(x_{i-1}(t),t)$}}
\newcommand{\uxxxi}{\mbox{$u_{xx}(x_i(t),t)$}}
\newcommand{\uxxxp}{\mbox{$u_{xx}(x_{i+1}(t),t)$}}
\newcommand{\uxxxm}{\mbox{$u_{xx}(x_{i-1}(t),t)$}}
\newcommand{\utxi}{\mbox{$u_t(x_i(t),t)$}}
\newcommand{\utxp}{\mbox{$u_t(x_{i+1}(t),t)$}}
\newcommand{\utxm}{\mbox{$u_t(x_{i-1}(t),t)$}}
\newcommand{\li}{\mbox{$l_i$}}
\newcommand{\lp}{\mbox{$l_{i+1}$}}
\newcommand{\lm}{\mbox{$l_{i-1}$}}
\newcommand{\suxs}{\mbox{$\sqrt{1\ +\ \ux\ }$}}
\newcommand{\suxsc}{\mbox{$\sqrt{(1\ +\ \ux)^3}$}}
\newcommand{\suxsf}{\mbox{$\sqrt{(1\ +\ \ux)^5}$}}
\newcommand{\ouxs}{\mbox{$\ 1\ +\ \ux\ $}}
\newcommand{\ft}{\mbox{$f(\theta)$}}
\newcommand{\fpt}{\mbox{$f^{\prime}(\theta)$}}
\newcommand{\fppt}{\mbox{$f^{\prime\prime}(\theta)$}}
\newcommand{\fpppt}{\mbox{$f^{\prime\prime\prime}(\theta)$}}
\newcommand{\cone}{\mbox{$c_1$}}
\newcommand{\ctwo}{\mbox{$c_2$}}
\newcommand{\cthree}{\mbox{$c_3$}}
\newcommand{\cfour}{\mbox{$c_4$}}
\newcommand{\xu}{\mbox{$u_x$}}
\newcommand{\suxis}{\mbox{$\sqrt{1\ +\ \uti\ }$}}
\newcommand{\sWpp}{\mbox{$\suxs\ \WD (u_x)$}}
\newcommand{\s}{\mbox{$\displaystyle \sum_i$}}
\newcommand{\sm}{\mbox{$\displaystyle \sum_{i=1}^{N-1}$}}
\newcommand{\smp}{\mbox{$\displaystyle \sum_{i=2}^{N-1}$}}
\newcommand{\inte}{\mbox{$\displaystyle \int_0^1$}}
\newcommand{\integ}{\mbox{$\displaystyle \int_{x_{i-1}}^{x_i}$}}
\newcommand{\sumint}{\mbox{$\displaystyle \s \int_{x_{i-1}}^{x_i}$}}
\newcommand{\q}{\mbox{$\displaystyle\frac{\ \up\ -\ \um\ }{2\,\li}\ r_i$}}
\newcommand{\newq}{\mbox{$\displaystyle\frac{\up\ -\ \um}{2\,\li}\ r_i$}}
\newcommand{\qq}{\mbox{$\displaystyle \frac{\tui}{\ \ \q\ \ }$}}
\newcommand{\dif}{\mbox{$u_{xx}\ -\ \q$}}
\newcommand{\uxx}{\mbox{$u_{xx}$}}
\newcommand{\uxxs}{\mbox{$u_{xx}^2$}}
\newcommand{\Wbi}{\mbox{$\Wb^\prime(\uxxi)$}}
\newcommand{\Wbp}{\mbox{$\Wb^\prime(\uxxp)$}}
\newcommand{\Wbm}{\mbox{$\Wb^\prime(\uxxm)$}}
\newcommand{\Wbpppz}{\mbox{$\Wb^{\prime\prime\prime}(\eta_i(x,t))$}}
\newcommand{\ut}{\mbox{$\sWpp\ \uxx$}}
\newcommand{\seg}{\mbox{$-\ \sm \uxxi\ [\ \tup\ -\ \tui\ ]$}}
\newcommand{\ter}{\mbox{$-\ \s \ui\ [\ \utxi\ -\ \utxm\ ]$}}
\newcommand{\iut}{\mbox{$\sqrt{1\ +\ u_x^2(x_i(t),t)}\ 
                \WD(u_x(x_i(t),t))\ u_{xx}(x_i(t),t)$}}
\newcommand{\mut}{\mbox{$\sqrt{1\ +\ u_x^2(x_{i-1}(t),t)}\ 
                \WD(u_x(x_{i-1}(t),t))\ u_{xx}(x_{i-1}(t),t)$}}
\newcommand{\quarto}{\mbox{$-\ \integ \sWpp\ \uxxs\ dx$}}
\newcommand{\dWi}{\mbox{$\dW[\ui]$}}
\newcommand{\WDi}{\mbox{$\WD[\ui]$}}
\newcommand{\WbDi}{\mbox{$\Wb^{\prime\prime}[\ui]$}}
\newcommand{\ite}{\mbox{$\vartheta_i$}}
\newcommand{\imte}{\mbox{$\vartheta_{i-1}$}}
\newcommand{\tei}{\mbox{$\theta_i$}}
\newcommand{\teip}{\mbox{$\theta_{i+1}$}}
\newcommand{\teim}{\mbox{$\theta_{i-1}$}}
\newcommand{\difi}{\mbox{$\up\ -\ \ui$}}
\newcommand{\difc}{\mbox{$\up\ -\ \um$}}
\newcommand{\difct}{\mbox{$\displaystyle \frac{\difc}{2}$}}
\newcommand{\deltai}{\mbox{$\Delta_i$}}
\newcommand{\difm}{\mbox{$\ui\ -\ \um$}}
\newcommand{\Wti}{\mbox{$W^{\prime\prime\prime}(\tei)$}}
\newcommand{\iWt}{\mbox{$W^{\prime\prime\prime}(\ite)$}}
\newcommand{\imWt}{\mbox{$W^{\prime\prime\prime}(\imte)$}}
\newcommand{\xiii}{\mbox{\Large{$\Xi$}}}
\newcommand{\xiiii}{\mbox{$\xiii_i$}}
\newcommand{\II}{\mbox{$I\! I$}}
\newcommand{\IIi}{\mbox{$I\! I_i$}}
\newcommand{\IIoi}{\mbox{$(I\! I_1)_i$}}
\newcommand{\IIti}{\mbox{$(I\! I_2)_i$}}
\newcommand{\IIo}{\mbox{$I\! I_1$}}
\newcommand{\IIt}{\mbox{$I\! I_2$}}
\newcommand{\III}{\mbox{$I\! I\! I$}}
\newcommand{\IIIi}{\mbox{$I\! I\! I_i$}}
\newcommand{\IIIoi}{\mbox{$(I\! I\! I_1)_i$}}
\newcommand{\IIIti}{\mbox{$(I\! I\! I_2)_i$}}
\newcommand{\IV}{\mbox{$I\! V$}}
\newcommand{\IVi}{\mbox{$I\! V_i$}}
\newcommand{\I}{\mbox{$I$}}
\newcommand{\AAA}{\mbox{$A$}}
\newcommand{\BBB}{\mbox{$B$}}
\newcommand{\CC}{\mbox{$C$}}
\newcommand{\DD}{\mbox{$D$}}
\newcommand{\mal}{\mbox{$\Wb^\prime(u_x(x,t))$}}
\newcommand{\ap}{\mbox{$a^\prime$}}
\newcommand{\umto}{\mbox{$(u_t^m)_1$}}
\newcommand{\uto}{\mbox{$u_t(0,t)$}}
\newcommand{\umtn}{\mbox{$(u_t^m)_N$}}
\newcommand{\utn}{\mbox{$u_t(1,t)$}}
\newcommand{\umo}{\mbox{$u^m(0,t)$}}
\newcommand{\umn}{\mbox{$u^m(1,t)$}}
\newcommand{\umxo}{\mbox{$(u_x^m)_1$}}
\newcommand{\uxo}{\mbox{$u_x(0,t)$}}
\newcommand{\umxn}{\mbox{$(u_x^m)_N$}}
\newcommand{\umxnn}{\mbox{$(u_x^m)_{N-1}$}}
\newcommand{\umxnp}{\mbox{$(u_x^m)_{N+1}$}}
\newcommand{\umtnn}{\mbox{$(u_t^m)_{N-1}$}}
\newcommand{\uxn}{\mbox{$u_x(1,t)$}}
\newcommand{\lnn}{\mbox{$l_N$}}
\newcommand{\lo}{\mbox{$l_1$}}
\newcommand{\umxoo}{\mbox{$(u_x^m)_0$}}
\newcommand{\deltao}{\mbox{$\Delta_1$}}
\newcommand{\suxos}{\mbox{$\sqrt{1\ +\ (u_x^m)_1^2\ }$}}
\newcommand{\suxns}{\mbox{$\sqrt{1\ +\ (u_x^m)_N^2\ }$}}
\newcommand{\numxt}{\mbox{$(u_x^m)_2$}}
\newcommand{\numtt}{\mbox{$(u_t^m)_2$}}
\newcommand{\Wone}{\mbox{$W[\umxo]$}}
\newcommand{\Wzero}{\mbox{$W[\umxoo]$}}
\newcommand{\Wtwo}{\mbox{$W[\numxt]$}}
\newcommand{\Wn}{\mbox{$W[\umxn]$}}
\newcommand{\Wnm}{\mbox{$W[\umxnn]$}}
\newcommand{\Wot}{\mbox{$\ \Wtwo\ -\ \Wone\ $}}
\newcommand{\Woto}{\mbox{$\ \Wone\ -\ \Wzero\ $}}
\newcommand{\Wnnm}{\mbox{$\ \Wn\ -\ \Wnm\ $}}
\newcommand{\ot}{\mbox{$\ \numtt\ -\ \umto\ $}}
\newcommand{\nnm}{\mbox{$\ \umtn\ -\ \umtnn\ $}}
\newcommand{\bdr}{\mbox{$\displaystyle \frac{\Wot}{\ \numxt\ -\ \umxo\ }$}}
\newcommand{\bdrm}{\mbox{$\displaystyle \frac{\Woto}{\ \umxo\ -\ \umxoo\ }$}}
\newcommand{\bdrp}{\mbox{$\displaystyle \frac{\Wnnm}{\ \umxn\ -\ \umxnn\ }$}}
\newcommand{\deno}{\mbox{$\ \numxt\ -\ \umxo\ $}}
\newcommand{\denn}{\mbox{$\ \umxn\ -\ \umxnn\ $}}
\newcommand{\bdo}{\mbox{$\displaystyle \frac{\ot}{\deno}$}}
\newcommand{\bdn}{\mbox{$\displaystyle \frac{\nnm}{\denn}$}}
\newcommand{\argt}{\mbox{$\displaystyle -\ \arctan \frac{1}{\xu}$}}
\newcommand{\argtb}{\mbox{$\left ( \argt \right )$}}
\newcommand{\argti}{\mbox{$\displaystyle -\ \arctan \frac{1}{\ui}$}}
\newcommand{\argtbi}{\mbox{$\left ( \argti \right )$}}
\newcommand{\bli}{\mbox{$L_i$}}
\newcommand{\tci}{\mbox{$\tilde{c}_i^0$}}
\newcommand{\tcp}{\mbox{$\tilde{c}_i^1$}}
\newcommand{\tcm}{\mbox{$\tilde{c}_i^{-1}$}}
\newcommand{\tcip}{\mbox{$\tilde{c}_{i+1}^0$}}
\newcommand{\tcpp}{\mbox{$\tilde{c}_{i+1}^1$}}
\newcommand{\tcmp}{\mbox{$\tilde{c}_{i+1}^{-1}$}}
\newcommand{\tcim}{\mbox{$\tilde{c}_{i-1}^0$}}
\newcommand{\tcpm}{\mbox{$\tilde{c}_{i-1}^1$}}
\newcommand{\tcmm}{\mbox{$\tilde{c}_{i-1}^{-1}$}}
\newcommand{\xidot}{\mbox{$\displaystyle 
        -\ \frac{\tup - \tui}{\up - \ui}$}}
\newcommand{\xmdot}{\mbox{$\displaystyle 
        -\ \frac{\tui - \tum}{\ui - \um}$}}
\newcommand{\xodot}{\mbox{$\displaystyle 
        -\ \frac{\numtt}{\numxt - \umxo}$}}
\newcommand{\xndot}{\mbox{$\displaystyle 
        \frac{\umtnn}{\umxn - \umxnn}$}}

\title[A crystalline algorithm for graphs moving by weighted curvature]{Convergence 
of a crystalline algorithm
for the heat equation in one dimension
and for the motion of a graph
by weighted curvature}

\author{Pedro Martins Gir\~{a}o}\thanks{Partially supported by a Dean's
Dissertation Fellowship from New York University and by
AFOSR grant 90-0090.} 

\address{Courant Institute,
        251 Mercer Street,
        New York, NY 10012}

\email{girao@acf9.nyu.edu} 

\curraddr{Mathematics Department,
Instituto Superior T\'{e}cnico,
1049-001 Lisbon, Portugal}

\email{pgirao@math.ist.utl.pt}

\author{Robert V. Kohn}\thanks{Partially supported by NSF
grant DMS-9102829, AFOSR grant 90-0090, and ARO contract DAAL03-92-G-0011.}

\address{Courant Institute,
        251 Mercer Street,
        New York, NY 10012}

\email{kohn@math5.nyu.edu}
\email{kohn@cims.nyu.edu}


\begin{abstract}
Motion by (weighted) mean curvature is a geometric evolution law for
surfaces, representing steepest descent with respect to (an)isotropic
surface energy.  It has been proposed that this motion could
be computed by solving the analogous evolution law using a
``crystalline'' approximation to the surface energy.  We present the
first convergence analysis for a numerical scheme of this type. Our
treatment is restricted to one dimensional surfaces (curves in the
plane) which are graphs. In this context, the scheme amounts to a new
algorithm for solving quasilinear parabolic equations in one space
dimension.
\end{abstract}

\maketitle

\section{Introduction}

In the modeling of phase transformations it is often of interest to
consider surface-energy-driven motion of interfaces (see
the work of Gurtin~\cite{G1} \cite{G2} \cite{G}, the recent review 
by \tch\ \cite{TCH}, and the references therein). 
If the surface energy is isotropic
this leads to motion by mean curvature, i.e.\ the normal velocity
of the interface equals its mean curvature. If the surface energy
is anisotropic the associated evolution equation has been called
``motion by weighted mean curvature'' (see the review Taylor~\cite{TW}).

        When the surface energy is ``strictly convex'' the evolution law is a
quasilinear parabolic equation. The isotropic case falls in this
class. When the surface energy is ``crystalline'' the surface must be
faceted and its evolution law reduces to a family of ordinary
differential equations for the lengths of the faces.  
The mathematical theory of surface-energy-driven motion of interfaces is
discussed at length in Angenent and Gurtin \cite{AG},
developing upon an extensive materials science literature of which
Herring \cite{H} is representative.  The essential aspects are
summarized in an appendix, for the reader's convenience.  We emphasize
that familiarity with this theory is {\em not} assumed in the present
work, though it provides the motivation and the context for
what we do.

        There has recently been intense mathematical activity 
concerning the analysis of motion by (weighted) mean curvature
(see for example 
Angenent and Gurtin~\cite{AG} \cite{A},
Brakke~\cite{Brakke},
Chen, Giga, and Goto~\cite{Chen},
de Mottoni and Schatzman~\cite{Mottoni},
Evans and Spruck~\cite{Evans} \cite{Spruck},
Gage and Hamilton~\cite{Gage},
Grayson~\cite{Grayson},
Gurtin~\cite{G},
Huisken~\cite{Huisken},
and Sethian~\cite{Sethian}).
The analysis of crystalline surface energies has also received considerable
attention both with regard to statics (energy minimization) (see
Sullivan \cite{S} and Taylor \cite{Ta1}) and
with regard to dynamics (motion by weighted mean curvature) (see 
Almgren, Taylor and Wang~\cite{ATW}, Angenent and Gurtin~\cite{AG}, 
Ohnuma and Sato~\cite{OS}, and Roosen and Taylor~\cite{RT}). 
We provide these references for the interested reader;
the only one we actually need in this \cpaper\ is~\cite{AG}.

It is a geometrically natural idea to approximate a strictly convex energy by
a crystalline one. 
This idea has recently been
analyzed by Sullivan in \cite{S}, but only for problems of energy minimization.
``Crystalline approximation'' has been used 
in applications without rigorous justification
as a method for computing 
motion by weighted mean curvature in the context 
of closed curves and surfaces, because
it replaces a parabolic differential equation
with a system of ordinary differential equations. 
This \cpaper\ represents a first attempt to establish its convergence.
We study the simplest possible case:
the interface is the graph of a function of one variable. Our
main result is the convergence of a ``crystalline'' scheme for solving
a quasilinear parabolic equation in one space
dimension. We prove convergence in $H^1$, with a specified rate. Our
method is somewhat similar to the convergence analysis for a Galerkin
approximation (see e.g. \thom\ \cite{T}).

        We emphasize that the ``crystalline'' discretization considered here
is quite different from any of the more standard schemes from numerical
analysis. The approximate solution is a piecewise linear function of
the spatial variable, with ``pieces'' that have fixed slopes and
variable lengths. This is in sharp contrast to finite difference and
finite element methods (where it is the spatial grid that is fixed),
to front tracking (which places no constraint
on the slopes of the approximate solution), and to spectral methods.

        The recent paper of Fukui and Giga \cite{FG} is closely
related to the present work.  They prove a general existence and
uniqueness theorem for motion of graphs by weighted curvature (by
adapting the theory of nonlinear semigroups), which applies to a wider
class of interfacial energies and even to ``incompatible'' initial
data.  When the surface energy is ``crystalline'' and the initial data
is piecewise linear and compatible (i.e.\ has the ``right'' slopes)
their solution is the same as ours.  They prove continuous dependence
on the initial data and on the form of the interfacial energy.  This
fact would suffice to prove convergence of our approximation scheme.
The result we prove here is sharper, however, because we get a
specific convergence rate in the $H^1$\/ norm.

        In another paper \cite{Girao} we prove convergence 
of the crystalline algorithm for the motion of a simple closed
convex curve by weighted curvature.  This is done by parametrizing
the weighted curvature by the angle between the interface normal
and a fixed coordinate axis, and comparing it with
the weighted curvature of the approximate solution.  It turns out that
the crystalline approximation scheme corresponds to a
standard finite difference scheme for the (nonlinear) evolution equation for
the weighted curvature.

        The organization of this \cpaper\ is as follows: in \cSection~2 we
set up the scheme and make some preliminary observations.  In
\cSection~3 we study the case of the ordinary heat equation. In
\cSection~4 we prove convergence for general convex energies and
constant Dirichlet boundary conditions.  In \cSection~5 we prove
convergence for constant Neumann boundary conditions.  In \cSection~6 we
set up the general Dirichlet problem.  However, for general Dirichlet
boundary conditions we do not prove convergence of the scheme because
we have been unable to prove that the approximate solution exists up
to a fixed time as the discretization gets finer.  It seems plausible
that, under suitable growth conditions on the energy, one
should be able to bound the slope of the approximate solutions.  Then
convergence would follow as in \cSection~4.  Finally, \cSection~7 is an
appendix which outlines the physical and mathematical context of our
work.

        This paper is part of the first author's Ph.D. thesis.

\section{Setup and preliminary observations}

        This paper 
is concerned with the convergence of an approximation scheme for the 
equation\addtocounter{equation}{+1}\setcounter{curvature}{\theequation}
$$ 
\left \{ \begin{array}{l}\displaystyle
        \frac{u_t}{\suxs}\ =\ 
                \WD(u_x)\ u_{xx}\qquad \mbox{if} \ \ 0\leq x \leq 1\\
                u(x,0)\ =\ u_0(x) 
        \end{array} \right. .
\eqno{(\thecurvature .1)}
$$
For now we focus on the case of homogeneous Dirichlet boundary
conditions,
$$
        u(0,t)\ =\ u(1,t)\ =\ 0 ,
\eqno{(\thecurvature .2)}
$$
postponing consideration of other boundary conditions to \cSections~5 and 6.
If $W(u_x)\ =\ \suxs$,\ \ $\WD(u_x)\ u_{xx}$ \ is the curvature of the graph
of $u$, and in general it is the 
negative of the gradient of 
\[ E(u)\ \stackrel{\triangle}{=}\ \inte W(u_x)\ dx ,\]
or the weighted curvature of the graph of $u$.
So Eq.~(\thecurvature) says that the normal velocity of the
graph of $u$ equals its weighted curvature.

        We assume throughout that $W$ is strictly convex and $C^3$, and that
Eq.~(\thecurvature) has a $C^3$ solution in $[0, 1]\times [0, T]$.
The existence of a $C^3$ solution is guaranteed by
Theorem~$\mbox{V}\!\mbox{I}.5.2$
and the regularity results of \lsu\ \cite{L} when $W$\/ and $u_0$\/ are 
sufficiently smooth, $u_0$\/ satisfies compatibility conditions at
zero and one, and $W$\/
satisfies certain growth conditions.  Appropriate growth conditions
are discussed in \cSection~6 (see especially conditions~(32)).

\vspace{\baselineskip}

The numerical scheme considered in 
this paper 
arises by
deriving an analogue of Eq.~(\thecurvature) in the setting where
$W$\/ is substituted by a piecewise linear function,
\Wt, which coincides with $W$\/ at its corners (see Figure~1).
We now carry out this derivation.  As explained in the Appendix,
this has been done in a slightly different way by Gibbs
(see the review Taylor~\cite{TW}).
The equations we will arrive at have been obtained
in a physical context by Angenent and Gurtin \cite{AG}.

\includegraphics{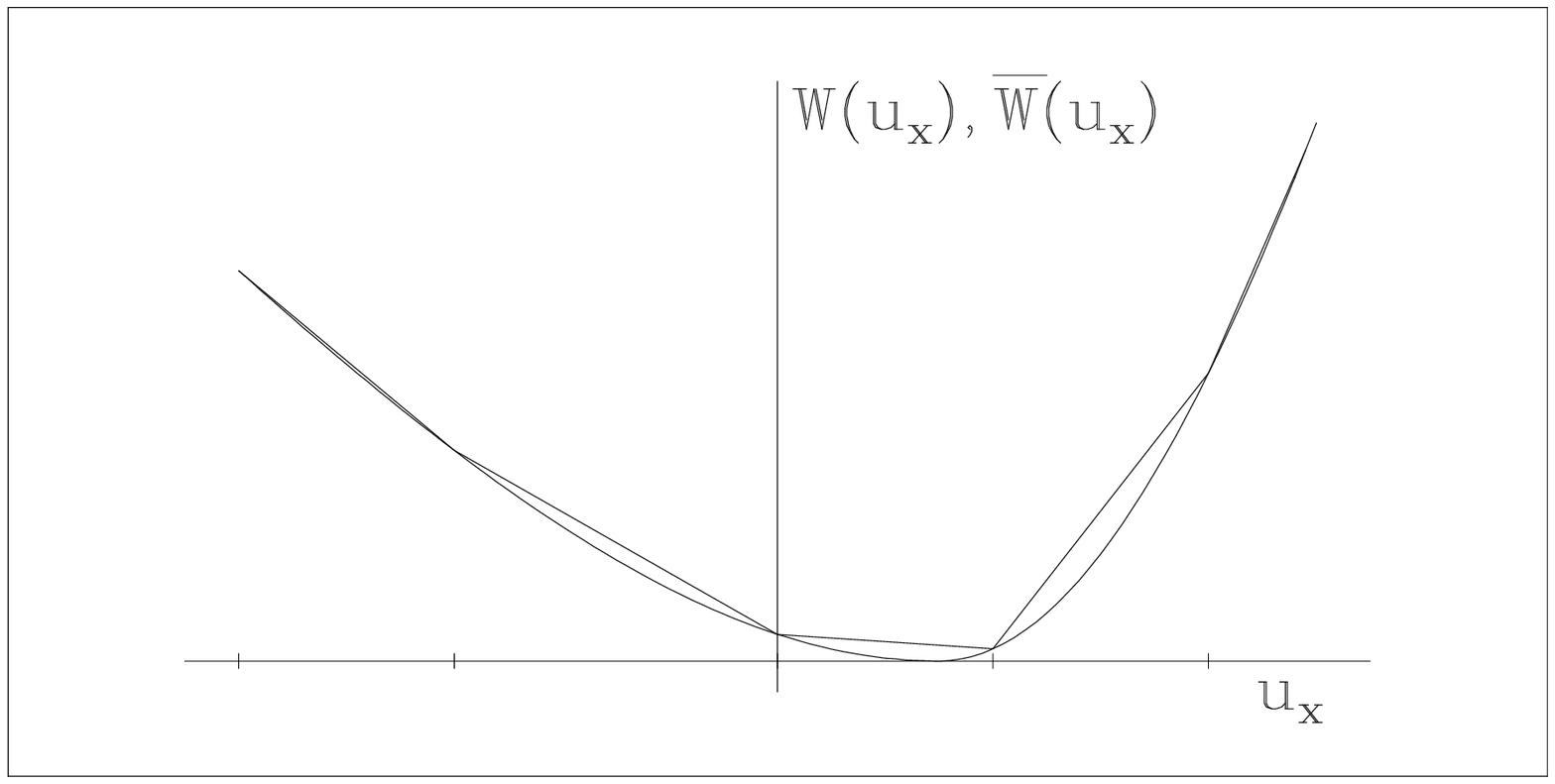}
\vglue 1.6truein                                              

\begin{center}
        Figure 1.\label{fone} The energy density $W$ approximated by a\\
        piecewise linear function \Wt. The function \Wt\\
        coincides with $W$ at its corners.
\end{center}

        The approximation scheme should deal with continuous 
piecewise linear functions,
$u^m$, where each segment, which from now on will be called face,
has an admissible slope, i.e.\ has slope equal to 
a value of $u_x$ corresponding to one of the corners of \Wt.
We will refer to the set of all such functions as the set of
admissible functions.
For now, the parameter $m$ has no meaning by itself, but later on it will
be the maximum distance between two such adjacent values of $u_x$,
and it will determine the rate of convergence of the scheme,
thus playing a role similar to the one that the mesh size does
in a finite difference scheme.

        We order the faces from left to right. 
At selected times faces will disappear and each time this happens we 
reorder them.  So the number
of faces is a piecewise constant function of time
denoted simply by $N$. 
We usually keep the time dependence of such functions
implicit.  The exception will be $u$'s time dependence and
that of its first argument, when that argument is itself a function,
which we often make explicit.
Let $\ x_i\ $ ($0\ =\ x_0\ <\ x_1\ \cdots\ <\ x_{i-1}\ <\ x_i\ 
<\ \cdots\ <\ x_{N-1}\ <\ x_N\ =\ 1$)
be the abscissa of the $i$th corner of $u^m$, 
$l_i$ be the length of the
projection of the $i$th face on the $x$-axis, and
\ui\ and \tui\ be the values of $u_x^m$ and $u_t^m$, 
respectively, for the $i$th face of $u^m$. 

        First of all, we focus on the initial data.     
Define $u_0^m(x)\ \stackrel{\triangle}{=}\ u^m(x,0)$.  We assume that

\vspace{2 mm}

\noindent
{\em $u_0^m$ belongs to the set of admissible functions, 
satisfies the boundary conditions, 
and jumps in $((u_0^m)_x)_i$ correspond to
adjacent corners of \Wt}\ (see Figure~2).\hfill
\addtocounter{equation}{+1}
\setcounter{keep}{\theequation}
(\arabic{equation})

\vspace{2 mm}

\noindent This is a reasonable assumption, 
and we make it because our scheme will move but not
create faces. 

\includegraphics{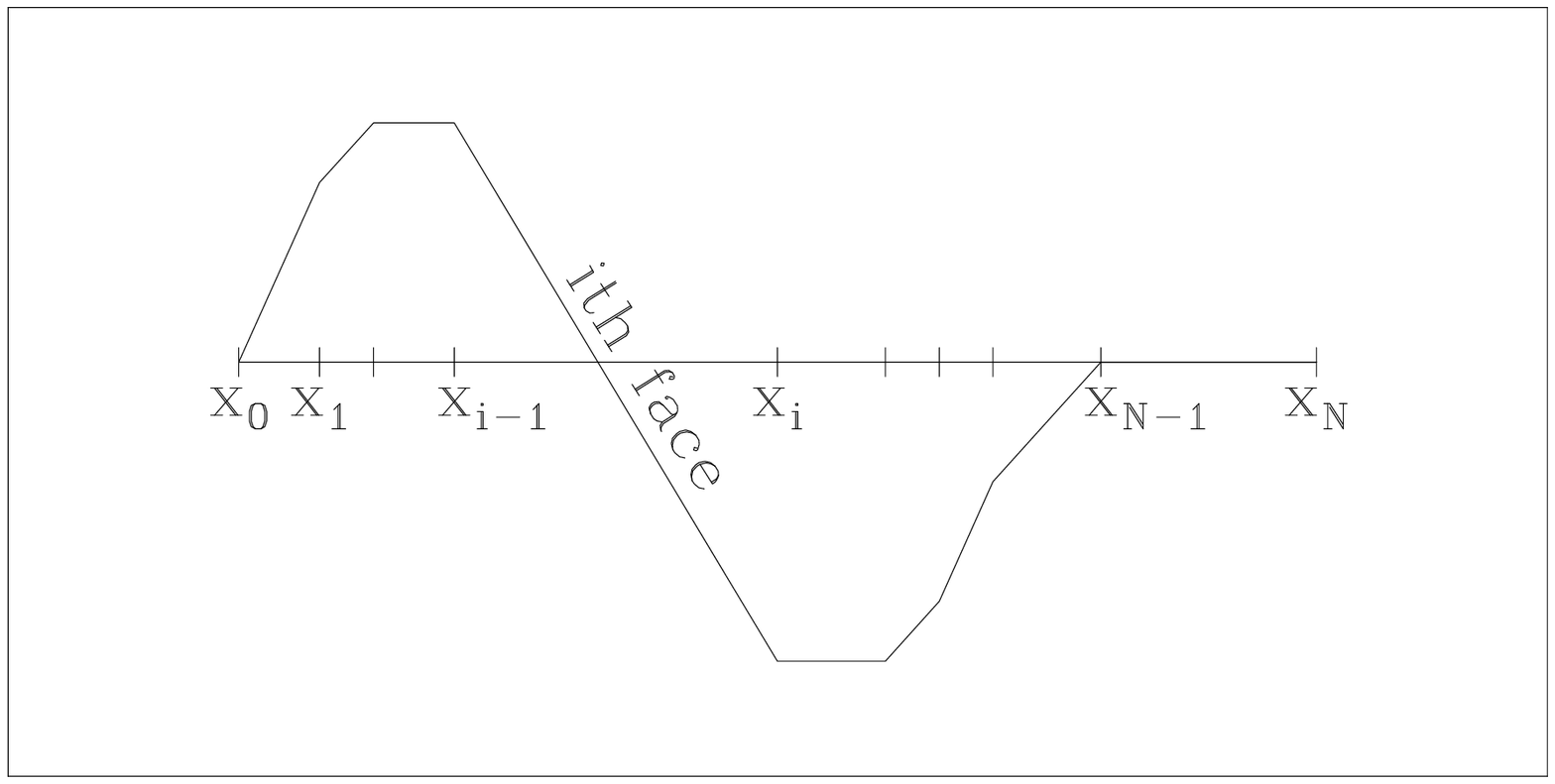}
\vglue 1.6truein                                              

\begin{center}
        Figure 2.\label{ftwo} A function $u_0^m(\ \cdot\ )$.
\end{center}

        Second, we impose the boundary conditions.  Since 
$u(0,t)\ =\ u(1,t)\ =\ 0$ we have to choose
\begin{equation}
(u_t^m)_1\ =\ 0\qquad \mbox{and}\qquad (u_t^m)_N\ =\ 0 \label{bdryut}
\end{equation} 
for them to be satisfied. 

        Our next goal is to derive the evolution equation of $u^m$.
Before doing this, however, we use the fact that $u^m$\/ is continuous
to derive some important relations.  The functions $u_x^m$\/
and $u_t^m$\/ are discontinuous, but the vertical velocity
of the $i$th corner of $u^m$ is the same whether computed
using the $i$th or $(i\ +\ 1)$th face. So,
differentiating $u^m(x_i(t),t)$\/ with respect to $t$,
\begin{equation}
        \xid\ =\ -\ \frac{\ \tup\ -\ \tui\ \ }{\ \up\ -\ \ui\ \ },
        \qquad\mbox{for}\ \ 1\ \leq\ i\ \leq\ N\ -\ 1 .
\label{xdot}
\end{equation}
It is sometimes more convenient to work with the lengths of the
faces rather than with the positions of the corners of $u^m$.  
Eqs.~(\ref{xdot})
imply \addtocounter{equation}{+1} \setcounter{ldot}{\theequation}
$$
\begin{array}{lcccl}
        \dot{l}_1\ &=\ & \dot{x}_1\       &=\ &-\ \bdo,\\ &&&&\\
        \dot{l}_N\ &=\ &-\ \dot{x}_{N-1}\ &=\ &\ \ \ \bdn,\end{array}
\eqno(\theldot .1^\prime)
$$
but these can be simplified using the boundary conditions.  We get
$$ \dot{l}_1\ =\ \dot{x}_1\ =\ c_1^1 (u_t^m)_2, \qquad \dot{l}_N\ =\ 
                -\ \dot{x}_{N-1}\ =\ c_N^{-1} (u_t^m)_{N-1},
\eqno(\theldot .1)
$$
and also (see Eqs.~(10.18) of \cite{AG})
$$
 \dot{l}_i\ =\ \dot{x}_i\ -\ \dot{x}_{i-1}\ =\  \ci \tui\ 
                                        +\ \cm \tum\ +\  \cp \tup,  
\eqno(\theldot .2)$$
for $2\ \leq\ i\ \leq\ N - 1$, with 
\begin{equation}\begin{array}{lcl}
        \ci\ & = &\displaystyle
                \frac{1}        
        {\ \ui\ -\ \um\ } \ 
                +\ \frac{1}
        {\ \up\ -\ \ui\ \ } , \\ &&\\
        \cm\ &=&\displaystyle
             -\ \frac{1}
        {\ \ui\ -\ \um\ } ,\\ &&\\ 
        \cp\ &=&\displaystyle -\ \frac{1}
        {\ \up\ -\ \ui\ \ } .\end{array}\label{KK}\end{equation}

        Finally, we turn to determining the evolution equation of $u^m$.
We compute the normal velocities of the faces of
$u^m$ in terms of the $l_i$'s.  
They should equal the negative of the gradient of 
\[
        \bar{E}(u)\ \stackrel{\triangle}{=}\ \inte \bar{W}(u_x)\ dx ,
\]
on the space of admissible functions, at $u^m$.  
By the definition of $\bar{W}$, 
\[ \Et(u^m)\ =\  \sum_{i=1}^N \Wi\ l_i . \]
We now compute the negative of the
gradient of \Et. 
Eqs.~(\ref{xdot}) imply that the rate of change of energy ``at'' 
the $i$th corner is
\[ -\ \{\ \Wi\ -\ \Wp\ \} \cdot \frac{\ \tup\ -\tui\ \ }{\ \up\ -\ui\ \ }\]
and, summing in $i$ from $1$ to $N-1$, that the (total) rate of change 
in energy is
\[ -\ \smp \ \left\{ \ \ppW\ \right\} \ \cdot\ \tui\ 
\]
\[ -\ \bdr\ \cdot\ \umto\ +\ \bdrp\ \cdot\ \umtn\ 
. \]
The boundary terms cancel because of the boundary conditions.

The preceding expression is the directional derivative of \Et,\ in the 
direction of the piecewise constant function which has value \tui\ 
for $x$ in $(x_{i-1}, x_i)$, at $u^m$.  In other words,
it is equal to
\[
        \left . \frac{d}{dh}\ \bar{E}(\phi_h)\ \right |_{h=0} ,
\]
where $\phi_h$\/ is admissible, $\phi_0\ =\ u^m$, and
\[
        \left .\frac{d\phi_h}{dh}\ \right |_{h=0}
                =\ \tui
\]
for $x$\/ in $(x_{i-1},x_i)$.  

Our space of admissible functions is somewhat unfamiliar,
so perhaps some elaboration is appropriate.
Suppose $\phi_h$\/ is a curve in the space of admissible
functions such that $\phi_0\ =\ u^m$.  
Since $\phi_h$\/ is admissible and $\phi_0\ =\ u^m$,
we see that $\left .\frac{d\phi_h}{dh}\ \right|_{h=0}$\/ is a piecewise
constant function which is constant, with value say $\tui$,
in the interval $(x_{i-1},x_i)$.  Conversely, given a piecewise
constant function which is constant, say with value $\tui$,
in the interval $(x_{i-1},x_i)$\/ there exists a family of
admissible functions $\phi_h$\/ such that $\phi_0\ =\ u^m$\/ and
$\left .\frac{d\phi_h}{dh}\ \right|_{h=0}\ =\ \tui$\/ in $(x_{i-1},x_i)$.
To obtain such a $\phi_h$\/ one can, for example, translate the $i$th
face of $u^m$\/ (for each $i$) vertically by $h\tui$.  It follows that
directional derivatives of $\bar{E}$\/ at $u^m$\/ on the space of admissible
functions are directional derivatives in the direction of
piecewise constant functions which are constant in $(x_{i-1},x_i)$.

        So the negative of the
gradient of \Et\ at $u^m$ is a piecewise constant function
that has value 
\[ \frac{\Delta_i}{l_i}, \qquad 
        \mbox{for $x$ in $(x_{i-1}, x_i)$},
\]
\begin{equation}
         \Delta_i \ \stackrel{\triangle}{=} \ \ppW\ ,
\label{delta}
\end{equation}
and for $2\ \leq\ i\ \leq\ N - 1$,
and has value zero for $x$ in $(x_0, x_1)$ and\ \ $(x_{N-1}, x_N)$.
We conclude that,
for $2\ \leq\ i\ \leq\ N - 1$, the $i$th face should have 
\begin{equation}
\tui\ \ =\  \suxis\ \cdot\ \frac{\Delta_i}{l_i} .
\label{K}\end{equation}

        In summary, {\em $u^m$\/ is determined by
solving\/ {\em Eqs.~(\theldot)}\/ with the $c_i^j$'s given by\/
{\em Eqs.~(\ref{KK})}, the 
$(u_t^m)_i$'s by\/ {\em Eqs.~(\ref{bdryut}) and~(\ref{K})}, 
and the $\Delta_i$'s
by\/ {\em Eqs.~(\ref{delta})}; the initial data should satisfy 
condition\/~{\em (\thekeep)}.}
This amounts to solving a nonlinear system  of ordinary
differential equations for the $l_i$'s.

        We make here a few comments.\label{comments}
\begin{itemize}
\item For reference we note the following elementary ``summation by
parts'' formulas:
{\tiny
\begin{eqnarray*}
   \sum_{i=1}^{N-1} a_i (b_{i+1}-b_i) &=& 
    -\ \sum_{i=2}^{N-1} b_i (a_i-a_{i-1})\ -\ a_1 b_1\ +\ a_{N-1} b_N ,\\
   \sum_{i=1}^N a_i (b_i-b_{i-1}) &=& 
    -\ \sum_{i=1}^{N-1} b_i (a_{i+1}-a_i)\ -\ a_1 b_0\ +\ a_N b_N ,\\
   \sum_{i=1}^{N-1} a_i (b_{i+1}-b_i) &=& 
    -\ \sum_{i=1}^N b_i (a_i-a_{i-1})\ -\ a_0 b_1\ +\ a_N b_N .
\end{eqnarray*}
}
\item Henceforth we write $\sum_i$ for $\sum_{i=1}^N$.
\item Later we will use the fact that if \up\ =\ \um\ then $\Delta_i\ =\ 0$.
\item If $W(\ \cdot\ )\ =\ \sqrt{\ 1\ +\ (\ \cdot\ )^2\ }$ then
$\frac{\Delta_i}{l_i}$ may be interpreted as the ``curvature of $u^m$
at its $i$th segment.''

\end{itemize}

\vspace{\baselineskip}
\vspace{\baselineskip}

        We now make some preliminary observations about this scheme.

\vspace{\baselineskip}

\noindent {\bf First},\label{first} note that one can 
think of $u^m$ as being defined on the whole real
line, the extension being odd with period two. 
Then\ $(u_x^m)_0\ =\ (u_x^m)_2$\ implies\ $(u_t^m)_1\ =\ 0$,\ and\ 
$(u_x^m)_{N+1}\ =\ (u_x^m)_{N-1}$\ implies\ $(u_t^m)_N\ =\ 0$.
This is important because it shows that for certain qualitative
features (such as our next observation) there is no need to
consider the end intervals separately.  It will also be useful
in our computations when we want to extend the range of some summations
by adding and subtracting boundary terms.  For example, in the expression
for the total rate of change of energy we can let the sum run
from 1 to $N$. Then the new boundary term at $x\ =\ 0$\/ is
\[ -\ \bdrm\ \cdot\ \umto. \]
Remark:\label{rmk}\
Eqs.~(\ref{delta}) and (\ref{K}) now hold for $i$ equal to 1 and $N$.

\vspace{\baselineskip}

\noindent {\bf Second}, we show that 

\vspace{2mm}

\noindent {\em even as time evolves and certain faces
disappear, jumps in \ui\ correspond to
adjacent corners of \Wt,\label{back}}
\hfill
\addtocounter{equation}{+1}
\setcounter{jump}{\theequation}
(\arabic{equation})

\vspace{2 mm}

\noindent and

\vspace{2mm}

\noindent {\em the maximum principle holds for $u_t^m$,}
\hfill
\addtocounter{equation}{+1}
\setcounter{maximum}{\theequation}
(\arabic{equation})

\vspace{2mm}

\noindent i.e.\ $u_t^m$\/ is bounded above by its maximum at time zero.
Statement~(\thejump) has also been noted by Taylor~\cite{TT}.
Notice the simple but important consequence of statement~(\thejump) 
that if \ui\ is smaller or greater than both \um\ and \up\ then\ \um\ =\ \up.\ 
Preparing for the proof, we note that 
if statement~(\thejump) is true at time $t$ then,
at that time, 
\[
\ci\ \tui\ \geq\ 0,\qquad \cm\ \tum\ \leq\ 0,\qquad \cp\ \tup\ \leq\ 0,
\]
since
$\Delta_i$, and consequently \tui, has the
same sign as $\up\ -\ \um$, due to the convexity of $W$. 
We proceed in four steps:
\begin{itemize}
\item {\em Statement\/~{\rm (\thejump)} holds at time zero, by condition\/~{\rm (\thekeep)}.
It clearly continues to hold until the time $t_1$\/
when one or more faces disappear.}
\item {\em The maximum principle
holds for $u_t^m$\/ for $t\ <\ t_1$.}\/
This is a consequence of the fact that
if $\max_{1\leq i\leq N}\/\tui$\/ is positive
then it is non-increasing.  For the proof, suppose the maximum
occurs at the $i$th face.  Then $\up\ -\um\ >\ 0$.  This
together with statement~(\thejump) implies
$\um\ <\ \ui\ <\ \up$; the first of these inequalities gives
$\cm\ <\ 0$\/ and the second $\cp\ <\ 0$.  It follows that
the right hand side of Eq.~($\theldot .2$) is
bounded below by $\ci\,\tui\ +\ \cm\,\tui\ +\ \cp\,\tui\ =\ 0$.
The left hand side of Eq.~($\theldot .2$) is 
\[
        \dot{l}_i\ =\ -\ \suxis\ \frac{\Delta_i}{\tuis}\ (u_t^m)_i^{\prime} ,
\]
by Eq.~(\ref{K}), where $(u_t^m)_i^{\prime}$\/ is the derivative
of $\tui$\/ ($(u_t^m)_i^{\prime}\ =\ (u_{tt}^m)_i$).  Therefore 
$(u_t^m)_i^{\prime}\ \leq\ 0$.  
\item {\em Statement\/~{\rm (\thejump)} holds at time $t_1$.}\/
Appealing again to Eq.~(\ref{K}), we see that only 
inflection faces can disappear, and the $l_i$'s have well defined
limits at time $t_1$.
Of course, no three or more adjacent
faces $i, \ldots, j$\/ can disappear at time $t_1$\/ for, if that were
the case, 
by Eqs.~($\theldot .2$), $\dot{l}_{i+1}\ =\ \cdots\ =\ \dot{l}_{j-1}\ =\ 0$\/
for $t\ <\ t_1$.  
At time $t_1$\/ two cases are possible:
\begin{itemize}
\item Case $(i)$\/ {\em The face $i$ disappears without its adjacent
neighbors disappearing.}  Then $\Delta_i\ =\ 0$\/ and
$\up\ =\ \um$.  At time $t_1$\/ the faces $i\ -\ 1$\/
and $i\ +\ 1$\/ join to from a single face.  (Note that 
in this case, for $t\ <\ t_1$,
$\Delta_{i-1} \Delta_{i+1}\ \leq\ 0$,
since either we have
$\um\ <\ \ui$\/ and $\ui\ >\ \up$,
or else the opposite of both of these inequalities.)
\item Case $(ii)$\/ {\em The faces $i\ -\ 1$ and $i$ disappear 
without their adjacent neighbors doing so.} 
Then $\Delta_{i-1}\ =\ \Delta_i\ =\ 0$, $(u_x^m)_{i-2}\ =\ \ui$\/ and   
$\um\ =\ \up$, and $(u_x^m)_{i-2}$\/ and\ \up\ correspond to adjacent
corners of \Wt.
\end{itemize}
\item {\em If the maximum of $u_t^m$\/ is discontinuous at time $t_1$\/
then it can only jump down.}\/  In fact, consider what happens
in each of the cases above:
\begin{itemize}
\item Case $(i)$\/ If the $(i\ -\ 1)$th and $(i\ +\ 1)$th faces
join to form an inflection face then for this new face $u_t^m\ =\ 0$.
Otherwise the new face has a value of $\Delta$\/ 
different from zero.  If $\Delta_{i-1} \Delta_{i+1}\ <\ 0$
for $t\ <\ t_1$\/ then the new face has a $\Delta$
equal to $\Delta_{i-1}$\/ or $\Delta_{i+1}$\/ 
(and the length of the new face at time $t_1$\/
is greater than the limit of the lengths of
both $l_{i-1}$\/ and $l_{i+1}$\/ as $t$\/ approaches $t_1$\/ from below).
If $\Delta_{i-1} \Delta_{i+1}\ =\ 0$ for $t\ <\ t_1$\/ 
then the new face might have a $\Delta$ different from
$\Delta_{i-1}$ and $\Delta_{i+1}$.  But then its $\Delta$
must be equal to the value of $\Delta_j$ for some $j$\/ such that
the $j$th face joins the $(i-1)$th and $(i+1)$th ones at time $t_1$\/
to form the new face.  (In this case the length of the new face at 
time $t_1$\/ is greater than the limit of the length of $l_j$\/
as $t$\/ approaches $t_1$.)
\item Case $(ii)$\/  The faces $i\ -\ 2$\/ and $i\ +\ 1$\/ do not 
join to form a new face.  
\end{itemize}
Therefore at time $t_1$\/ the maximum of $u_t^m$\/ can only decrease.  
\end{itemize}
By applying the above reasoning repeatedly, we conclude that
statement~(\thejump) and the maximum principle for $u_t^m$\/
hold for all time.  The above shows in particular that $u^m$\/ exists
for all time.

\vspace{\baselineskip}

\noindent {\bf Third}, we note the discrete analogue of
$\frac{d}{dt}\ \frac{1}{2}\ \int_0^1\ u_x^2\ dx\ =\ 
-\ \int_0^1\ u_t\ u_{xx}\ dx\
+\ u_t(1,t)\ u_x(1,t)\ -\ u_t(0,t)\ u_x(0,t)\ $.
This is an identity valid for any function $u$.  The discrete version
applies to any admissible function $u^m$, i.e.\ one which is piecewise
linear in space with slopes $\ui$, such that the jumps in $\ui$\/
correspond to adjacent corners of $\bar{W}$.  The proof uses only
Eqs.~($\theldot.1^{\prime}$) and ($\theldot.2$), not the discretized
differential equation~(\ref{K}) or the Dirichlet boundary 
condition~($\theldot.1$):
\begin{eqnarray}
\nonumber
        \ddt\ \half\ \s \uti\ l_i\! \!\!  &=
        &\!\!\!
         \half\ \smp \uti\ [\ \ci \tui\ +\  \cm \tum\ +\  \cp \tup\ ]\\
\nonumber &&\\ 
\nonumber
        &&\qquad\qquad +\ \half\ (u_x^m)_1^2 c_1^1 (u_t^m)_2\ 
        +\ \half\ (u_x^m)_N^2 c_N^{-1} (u_t^m)_{N-1}\\ 
\nonumber &&\\ 
\nonumber
        &&\qquad\qquad -\ \half\ (u_x^m)_1^2 c_1^1 (u_t^m)_1\
        -\ \half\ (u_x^m)_N^2 c_N^{-1} (u_t^m)_N\\ 
\nonumber &&\\
\nonumber
        &=&\!\!\! \half\ \smp \tui\ [\ \ci \uti\ 
        +\  c_{i+1}^{-1} \utp\ +\  c_{i-1}^1 \utm\ ]\\
\nonumber &&\\ 
\nonumber
        &&\qquad\qquad -\ \half\ (u_x^m)_1^2 c_1^1 (u_t^m)_1\
        -\ \half\ (u_x^m)_N^2 c_N^{-1} (u_t^m)_N\\ 
\nonumber &&\\
\nonumber
        &&\qquad\qquad +\ \half\ (u_x^m)_2^2 c_2^{-1} (u_t^m)_1\
        \! +\ \half\ (u_x^m)_{N-1}^2 c_{N-1}^1 (u_t^m)_N\\ 
\nonumber &&\\
\label{Wdot}
        &=&\!\!\! -\ \half\ \smp \tui\ [\ \up\ -\ \um\ ]\\
\nonumber &&\\
\nonumber
        &&\qquad\qquad 
        -\ \half\ (u_t^m)_1\ [\ (u_x^m)_1\ +\ (u_x^m)_2\ ] \\
\nonumber &&\\
\nonumber
        &&\qquad\qquad
        +\ \half\ (u_t^m)_N\ [\ (u_x^m)_{N-1}\ +\ (u_x^m)_N\ ] .
\end{eqnarray}
For the solution of Eq.~(\thecurvature) we have
$\frac{d}{dt}\ \int_0^1\ u_x^2\ dx\ \leq\ 0$.
Similarly, for the solution of the discretized problem with
homogeneous Dirichlet boundary condition we see from Eq.~(\ref{Wdot})
that $\frac{d}{dt}\ \sum_i\ \uti\ l_i\ \leq\ 0$.

\vspace{\baselineskip}

\noindent {\bf Fourth}, $\bar{E}(u^m)(\ \cdot\ )$\/ is a continuous 
function of $t$\/ but is only piecewise $C^1$\/ because faces disappear
at selected times.  We will not mention this restriction again,
although it will come up repeatedly, since it is of no consequence 
for our analysis.

\vspace{\baselineskip}

\noindent {\bf Fifth} and last, the maximum principle holds for $u^m$ because 
an interior maximum will occur at the $i$th corner 
($1\ \leq\ i\ \leq\ N\ -\ 1$)
only if $\ui\ >\ \up$\/ and in this case both \tui\ and \tup\ are nonpositive.
\setcounter{transfer}{\theequation}

\section{Convergence for the heat equation}

\setcounter{equation}{\thetransfer}
Before proving convergence of the approximation scheme
for Eq.~(\thecurvature), let us treat the simpler case
of the corresponding approximation scheme for the heat equation,
\begin{equation} \left \{ \begin{array}{l}
        \mbox{\normalsize{$u_t$}}\ =\ 
                u_{xx}\qquad \mbox{if} \ \ 0\leq x \leq 1\\
                u(x,0)\ =\ u_0(x) \\
        u(0,t)\ =\ u(1,t)\ =\ 0 \end{array} \right. .
\label{heat}\end{equation}
In this case
\[ W(u_x)\ =\ \half\ \ux\ \]
and
\begin{equation} 
        \tui\ =\ \frac{\Delta_i}{\li}, \qquad 
                        \Delta_i\ =\ \frac{\upmum}{2} .
\label{nosquare}
\end{equation}
We shall estimate the growth of the $H^1$ norm of $\ u\ -\ u^m\ $.  

        The calculation in this \csection\ would be simplest if
we took full advantage of the boundary conditions, dropping
terms such as $\uto$\/ or $\umto$\/ which vanish due to 
Eqs.~$(\thecurvature .2)$\/ and (\ref{bdryut}).  This would
be inefficient because it would force us to repeat many
calculations in \cSections~5 and 6, where we consider other
boundary conditions.  We therefore proceed as follows:
\begin{itemize}
\item The Dirichlet boundary conditions $u(1,t)\ =\ 0$\/
and $\umtn\ =\ 0$\/ will be used at 
$x\ =\ 1$, and we shall drop terms which vanish
as a result.
\item The Dirichlet boundary conditions will also be used at 
$x\ =\ 0$, but we shall label ((i)~through~(vii))
the associated boundary terms the first time they appear for future
reference.  Then these terms will be dropped since they vanish due to
$u(0,t)\ =\ 0$\/ and $\umto\ =\ 0$.
\item In \cSections~5 and 6 we will deal with boundary conditions
other than the homogeneous one $u(0,t)\ =\ \umto\ =\ 0$.  Specialized
to the case of the heat equation, they lead to relations of the form
\begin{equation}
        (u_t^m)_1\ =\ \frac{\Delta_1}{l_1}\ r_1 , 
\label{zero}
\end{equation}
where $r_1$\/ is a function of time, $0\ <\ r_1\ \leq\ 1$,
and $\Delta_1$\/ is defined by Eq.~(\ref{delta}) for $i\ =\ 1$.
Notice that $\Delta_1$\/ depends on $\umxoo$, $\umxo$, and $\numxt$.
\item To impose the Dirichlet boundary condition $u^m(0,t)\ =\ 0$,
one should take $\umxoo\ =\ \numxt$, corresponding to the odd,
periodic extension of $u^m$.  This yields $\Delta_1\ =\ 0$\/
and $\umto\ =\ 0$, as noted above.  In this case we may as well set 
$r_1\ =\ 1$\/ since both sides of Eq.~(\ref{zero}) vanish.
\item It is convenient to unify Eqs.~(\ref{nosquare}) and (\ref{zero})
by writing
\[
        (u_t^m)_i\ =\ \frac{\Delta_i}{l_i}\ r_i , 
\]
for $1\ \leq\ i\ \leq\ N$, with $r_2\ =\ \cdots\ =\ r_N\ =\ 1$,
and $0\ <\ r_1\ \leq\ 1$.  Zero boundary conditions at $x\ =\ 1$\/
imply $\Delta_N\ =\ 0$, just as zero boundary conditions at $x\ =\ 0$\/
imply $\Delta_1\ =\ 0$.
\item At risk of redundancy, we emphasize:  for homogeneous Dirichlet
boundary conditions it suffices to set $u_t(0,t)\ =\ 0$, $\umto\ =\ 0$,
$r_1\ =\ r_2\ =\ \cdots\ =\ r_N\ =\ 1$, and $\umxoo\ =\ \numxt$ in 
the calculation that follows, and the boundary terms
(i)~through~(vii) all vanish.
\end{itemize}

        One further remark.  The right hand side of Eq.~(\ref{Wdot})
can be written
\begin{equation}\begin{array}{l}
        -\ \half\ \smp \tui\ [\ \up\ -\ \um\ ]\\
\\
        -\ \half\ \umto\ [\ \numxt\ -\ \umxoo\ ]\ r_1\\
\\
        -\ \half\ \umto\ [\umxoo\ r_1\ +\ \umxo\ +\ \numxt\ (1\ -\ r_1)].
\end{array}\label{Wdots}\end{equation}
For the case of the heat equation, discretized as above, this becomes
\[\begin{array}{l}
        -\ \smp \tuis\ l_i\ -\ (u_t^m)_1^2\ l_1\\
\\
        -\ \half\ \umto\ [\umxoo\ r_1\ +\ \umxo\ +\ \numxt\ (1\ -\ r_1)].
\end{array}\]

        Now we begin the convergence analysis.  We have
\[ \begin{array}{lcr}
        \ddt\ \half\ \inte |\ u_x\ -\ u_x^m\ |^2 \ dx\ &=& 
        \ddt\ \left [\ \half\ \inte \ux\ dx\ +\ 
        \half\ \s \uti\ \li\right. \\ &&\\
        &&\left. 
-\ \sumint u_x\ u_x^m\ \right ] \\ 
\end{array} \]
\[ \begin{array}{rl}
        =&-\ \inte u_t^2\ dx\ \underbrace{-\ \uto\ \uxo}_{\mbox{(i)}}\
        -\ \s \tuis\ \li\ \\ &\\
        &\underbrace{-\ \half\ \umto\ 
        [\umxoo\ r_1\ +\ \umxo\ +\ \numxt\ (1\ -\ r_1)]}_{\mbox{(ii)}}
\\&\\
        &-\ \ddt\ \s \ui\ [\ \uxi\ -\ \uxm\ ]\\&\\
        =&\displaystyle -\ \inte |\ u_t\ -\ u_t^m\ |^2\ dx\ 
                \underbrace{-\ 2\ \ 
                \inte u_t\ u_t^m\ dx}_{\AAA} \\&\\
         \qquad\qquad\qquad\qquad
  &\underbrace{-\ \s \ui\ [\ \uxxi\ \xid\ -\ \uxxm\ \xmd\ ]}_{\BBB}\\
 &\\
        &   \underbrace{-\ \s \ui\ [\ \utxi\ -\ \utxm\ ]}_{\CC} .
\end{array} \]
We have dropped the terms which vanish due to the homogeneous
boundary condition, as we announced we would do.
We rewrite $A$ and $B$:
\begin{eqnarray}
        A&=&\displaystyle -\ 2\ \s \tui\ [\ \uxxi\ -\ \uxxm\ ] ,
\nonumber\\
\nonumber&&\\
        B&=&\sm \uxxi\ [\ \up\ -\ \ui\ ]\ \xid\ 
        +\ u_x(0,t)\ (u_x^m)_1\ \dot{x}_0 
\nonumber\\
\nonumber&&\\
         &=&-\ \sm \uxxi\ [\ \tup\ -\ \tui\ ]\ 
        \nonumber\\&&\nonumber\\
        &=&\s \tui\ [\ \uxxi\ -\ \uxxm\ ]\ 
        \underbrace{+\ \umto\ \uxo}_{\mbox{(iii)}}.
\label{BB}
\end{eqnarray}
Note that we have used Eqs.~(\ref{xdot}).
So,
\[ \begin{array}{lcl}
        A\ +\ B &=&-\ \s \tui\ [\ \uxxi\ -\ \uxxm\ ]\ 
\end{array}\]
\[\begin{array}{cl}
        =&\s\ \left [\ \q\ \right ]\ 
                [\ -\ \uxxi\ +\ \uxxm\ ]\ \\&\\
        =&\displaystyle\half\ \s\ [\ \up\ -\ \ui\ ]\ 
                \left [ \ -\ \frac{\uxxi\ -\ \uxxm\ }{\li\ }\ \right ]\ r_i
                 \\ &\\ 
        &\displaystyle\!\!\!\!\!\!\!+\ \half\ \s\ [\ \ui\ -\ \um\ ]\ 
                \left [ \ -\ \frac{\uxxi\ -\ \uxxm\ }{\li\ }\ \right ]\ r_i .
\end{array}
\]
On the other hand,
\begin{eqnarray}
\label{eqc}     C&=&\sm \utxi\ \cdot \ [\ \up\ -\ \ui\ ]\ 
        \underbrace{+\ \uto\ \umxo}_{\mbox{(iv)}}\\ \nonumber&&\\ \nonumber
        &=&\half\ \s\ [\ \up\ -\ \ui\ ]\ \cdot\
                \left [\ \utxi\ \right ]
                \\ \nonumber && \\ \nonumber 
        &&\!\!\!\!\!\!\!+\ \half\ \s\ 
                        [\ \ui\ -\ \um\ ]\ \cdot\
                \left [\ \utxm\ \right ]
                \\ \nonumber && 
        \\ \nonumber&&\!\!\!\!\!\!\!\underbrace{-\ \half\ 
                \uto\ [\ \umxo\ -\ \umxoo\ ]}_{\mbox{(v)}} .
\end{eqnarray}
Hence,\ $A\ +\ B\ +\ C$\/ equals 
\[ \begin{array}{l}
        +\ \s\ \onefourth\ [\ \up\ -\ \ui\ ]\ \cdot \ 
                \left [\ u_{xxx}(\alpha_i(t),t)\ \right ]\ r_i\ l_i
\\ \\
        -\ \s\ \onefourth\ 
                [\ \ui\ -\ \um\ ]\ \cdot \ 
                \left [\ u_{xxx}(\beta_i(t),t)\ \right ]\ r_i\ l_i
\\ \\ 
        \underbrace{+\ \half\ \uto\ [\ \umxo\ -\ \umxoo\ ]\ 
                (1\ -\ r_1)}_{\mbox{(vi)}}
\\ \\
        +\ \half\ u_t(x_1(t),t)\ [\ \numxt\ -\ \umxo\ ]\ (1\ -\ r_1),
\end{array} \]
for some\ $\alpha_i\ \mbox{and}\ \beta_i\ $\
belonging to $(x_{i-1},\ x_i)$.  The last term in this sum
can be written as
\[ 
        \half\ u_t(x_1(t),t)\ [\ \numxt\ -\ \umxo\ ]\ (1\ -\ r_1)
\qquad\qquad\qquad\qquad\qquad\qquad\qquad
\] 
\vspace{1 mm}
\[\begin{array}{ccl}
\qquad\qquad\qquad\ \ &=&\underbrace{+\ 
        \half\ \uto\ [\ \numxt\ -\ \umxo\ ]\ 
                (1\ -\ r_1)}_{\mbox{(vii)}}\\ &&\\
        &&+\ \half\ [\ \numxt\ -\ \umxo\ ]\ u_{xxx}(\gamma_1(t),t)\ 
                (1\ -\ r_1)\ l_1 ,\\&&
\end{array}\]
for some $\gamma_1$ belonging to $(0, x_1)$. 
For zero boundary conditions each of (i) through (vii) vanishes.
We stop to note that our calculation was exact up to here.
We will now bound $A\ +\ B\ +\ C$.  The final bound for
the $L^2$\/ norm of $u_x\ -\ u_x^m$\/ will be sharp inasmuch as 
the next bound is optimum.
The terms corresponding to $\alpha_1$\/ and $\gamma_1$\/ combine
and we obtain
\[
        |\ A\ +\ B\ +\ C\ |\ \leq\ 
\half\ ||\ u_{xxx}\ ||_\infty (t)\ \max_{1\leq i\leq N + 1}\ 
        |\ \ui\ -\ \um\ |;
\]
here $||\ ||_\infty$\ is the supremum norm in $x$ on the interval $[0, 1]$.
Thus the accuracy of the discretization is governed by 
\begin{equation}
        m\ \stackrel{\triangle}{=} \ \max_{1\leq i\leq N + 1}\ 
        |\ \ui\ -\ \um\ |. \label{MM}
\end{equation}
The parameter $m$\/ controls the convergence of the scheme in much the 
way the mesh size controls the behavior of a finite difference approximation.
Notice that our estimates do not require any uniformity on the
quantities $|\ \ui\ -\ \um\ |$\/ as $i$\/ varies.

        To recapitulate, we have shown that if $u$\/ solves the
heat equation~(\ref{heat}) and $u^m$\/ solves its discretized
version, both with homogeneous Dirichlet boundary conditions, then
\[ \begin{array}{lcl}
        \ddt\ \half\ \! \inte\! |\ u_x\ -\ u_x^m\ |^2 \ \! dx\! &\!
        \leq &\! -\ \inte\!\! |\ u_t\ -\ u_t^m\ |^2\ \! dx\ 
                +\ \half\ m\ ||\ u_{xxx}\ ||_\infty (t),
\end{array} \]
Suppose that the discrete initial data $u_0^m$\/
have been chosen so that
\begin{equation}
        ||\ (u_0^m)_x\ -\ (u_0)_x\ ||_\infty\ \leq\ m .
\label{comeco}
\end{equation}
Then the error at time $t$\/ is
\begin{equation}
        \half\ \inte |\ u_x\ -\ u_x^m\ |^2 \ dx\
\leq\ \half\ m^2\ +\ \half\ m\ \int_0^t ||\ u_{xxx}\ ||_\infty (\tau)\ d\tau .
\label{moment}
\end{equation}
By Poincar\'{e}'s inequality the $L^2$\/ norm of $u\ -\ u^m$\/
is controlled by the $L^2$\/ norm of its spatial derivative.  We
conclude that
\[
        \sup_{0\leq t\leq T} ||\ u\ -\ u^m\ ||_{H^1([0,1])}\ \leq\ c\
                                        m^{\frac{1}{2}} ,
\]
with $c$\/ depending on $u$\/ and $T$\/ but not on $m$.

        Since we are assuming $u_0$\/ is $C^3$,\label{constroi}
one can choose $u_0^m$\/ such that
the hypothesis~(\ref{comeco}) is satisfied.  For example, one can
proceed as follows.  First, identify all the points where
$(u_0)_x$\/ is admissible and draw the tangents to $u_0$\/
through those points.  Next, draw a line through each inflection
point $\check{x}$\/ of $u^m$\/ with the admissible slope $\ui$\/
closest to the value of $(u_0)_x(\check{x})$\/ and such that 
$|\ \ui\ |\ >\ |\ (u_0)_x(\check{x})\ |$.  Then, draw a
line through $(0,0)$\/ with the admissible slope closest to
the value of $(u_0)_x(0)$\/ and greater than (respectively, less than)
$(u_0)_x(0)$\/ 
if $(u_0)_{xx}(0)\ <\ 0$\/ ($(u_0)_{xx}(0)\ >\ 0$); 
if zero is an inflection point of $u_0$\/
disregard this step.  Do a similar construction at the point one.  Finally,
obtain $u_0^m$\/ by the union of segments on the lines above.

        It is worth noting that, for relatively large times $t$,
error estimate~(\ref{moment}) is dominated by the term\ \,
$m \int_0^t || u_{xxx} ||_{\infty}(\tau)\ d\tau$, not by the
approximation of the initial data.  A cruder approximation
of the initial data, $|| (u_0^m)_x\ -\ (u_0)_x ||_{L^2([0,1])}^2\ =\ 
\mbox{O}(m)$, would not affect the order of convergence.

\vspace{2mm}

        {\bf Remark}.\label{boundary} We collect the boundary terms which 
entered in the calculation above because we will need them later.
They were:
\[\begin{array}{rcl}
        \mbox{(i)}&=&-\ \uto\ \uxo ,\\
&&\\
        \mbox{(ii)}&=&-\half 
                \umto\ [\umxoo r_1 + \umxo + \numxt (1 - r_1)],\\
&&\\
        \mbox{(iii)}&=& \umto\ \uxo ,\\
&&\\
        \mbox{(iv)}+\mbox{(v)}+\mbox{(vi)}+\mbox{(vii)}&=&\half
                \uto\ \ [\umxoo r_1 + \umxo + \numxt (1 - r_1)].
\end{array}\]
If\ \ $\umto\ =\ \uto$, as in \cSection~6, then
the sum (i)~through~(vii) is zero.  Alternatively, if $\umxo\ =\ \uxo$,
as in \cSection~5, then
\[\begin{array}{rcl}
        \mbox{(ii)}+\mbox{(iii)}\!\!&=&\!\!\half
        \umto [- \umxoo r_1 + \umxo - \numxt (1 - r_1)],\\
&&\\
\!\!\mbox{(i)}+\mbox{(iv)}+\mbox{(v)}+\mbox{(vi)}+\mbox{(vii)}\!\!&=&\!\!
                \half \uto\ [\umxoo r_1 - \umxo + \numxt (1 - r_1)].
\end{array}\]

\vspace{2mm}

        Before proceeding to the next \csection\ we make a final
remark.  If the distance between any two adjacent admissible
slopes is a fixed constant, i.e.\ the mesh in the $u_x$-axis of Figure~1
is uniform, then it might be possible to improve the 
bound for $A\ +\ B\ +\ C$.
In fact, if the $i$th is not an inflection face of $u^m$\/
then either $\um\ <\ \ui\ <\ \up$\/ or both these inequalities are
reversed.  In either case the absolute value of the $i$th term in the 
expression for the sum of $A\ +\ B\ +\ C$\/ is equal to
\[
        \onefourth\ m\ \left |\ u_{xxxx}(\mu_i(t),t)\ l_i\ 
                (\alpha_i(t)\ -\ \beta_i(t))\ \right | ,
\]
for some $\mu_i$\/ belonging to $(x_{i-1},x_i)$, so less than
\[
        \onefourth\ m\ ||\ u_{xxxx}\ ||_{\infty}(t)\ l_i^2 .
\]
If the $i$th is an inflection face of $u^m$\/ then the $i$th
term of $A\ +\ B\ +\ C$\/ is bounded by
\[
\half\ m\ ||\ u_{xxx}\ ||_{\infty}(t)\ l_i .
\]
Hence 
\[\begin{array}{lcl}
        |\ A\ +\ B\ +\ C\ |&\leq&\displaystyle \phantom{+}\ \onefourth\ m\ 
                ||\ u_{xxxx}\ ||_{\infty}(t)\ \s\ l_i^2 \\
&&\\
                &&\displaystyle +\ \half\ m\  ||\ u_{xxx}\ ||_{\infty}(t)
                \sum_{\stackrel{i\mbox{\tiny th an}}{\mbox{\tiny infl. face}}}
                l_i .
\end{array}\]
This suggests that to get a sharper estimate than the one obtained
above one should bound $\max_i l_i$\/ on an interval $[0,T]$\/ in terms 
of $m$.  Note that for the present case Eqs.~($\theldot .2$) reduce to
\[\begin{array}{lcll}
\left [ \begin{array}{c} \dot{l}_i\\ \dot{l}_{i+1}\\ \vdots \\ 
        \dot{l}_{j-1}\\ \dot{l}_j 
        \end{array} \right ]
&=&
\left [ \begin{array}{cccccc} 
        2       &-1     &0      &\cdots &0      &0\\ 
        -1      &2      &-1     &\cdots &0      &0\\ 
        \vdots  &\vdots &\vdots &\vdots &\vdots &\vdots\\ 
        0       &0      &0      &\cdots &2      &-1\\
        0       &0      &0      &\cdots &-1     &2
        \end{array} \right ]&
\left [ \begin{array}{c} \frac{1}{l_i}\\ \frac{1}{l_{i+1}}\\ \vdots \\ 
                \frac{1}{l_{j-1}}\\ \frac{1}{l_j} \end{array} \right ]
\end{array}
\]
between inflection faces.
\setcounter{transfer}{\theequation}

\section{Convergence for general convex energies}

\setcounter{equation}{\thetransfer}
We turn now to 
proving convergence of the approximation scheme described in \cSection~2,
for Eq.~(\thecurvature).  The argument is in many respects similar
to that for the heat equation.
Of course, the functions $u$\/ and $u^m$\/ 
are defined differently than they were in \cSection~3, because
we are solving different differential equations.  However,
we shall keep the notation and labels used in the previous
\csection\ to highlight the parallels between the two calculations.
We shall also keep the conventions concerning the handling of
boundary terms.  

        The main idea is, once again, to control
the evolution of the $H^1$ norm of $u\ -\ u^m$.  Arguing as
before, we get
\[ \begin{array}{lcr}
        \ddt\ \half\ \inte |\ u_x\ -\ u_x^m\ |^2 \ dx\ &=& 
        \ddt\ \left [\ \half\ \inte \ux\ dx\ +\ 
        \half\ \s \uti\ \li\right. \\ &&\\
        &&\left. 
        -\ \sumint u_x\ u_x^m\ \right ] \\ 
\end{array} \]
\begin{equation}\begin{array}{rl}
        =&-\ \inte \sWpp\ \uxxs\ dx\ +\ \mbox{(i)}\\
        &\\ &-\ \half\ \s \tui\ [\ \up\ -\ \um\ ]\ r_i\ +\ \mbox{(ii)}\\&\\
        &\underbrace{+\ \s \tui\ [\ \uxxi\ -\ \uxxm\ ]}_{\BBB}\ +\ \mbox{(iii)}
\\&\\
        &\underbrace{+\sm \utxi\ [\ \up\ -\ \ui\ ]}_{\CC}\ +\ \mbox{(iv)} . \\
\end{array}\label{hnorm} \end{equation}
We have used expression~(\ref{Wdots}) and Eqs.~(\ref{BB}) and (\ref{eqc}).
It is still convenient to rewrite $B$:
\[ 
        B\ =\ \sumint \tui\ \uxx\ dx . 
\]
It is easy to check that expression~(\ref{hnorm}) can be written as
the sum
\[ \begin{array}{l}
\underbrace{
-\sumint\! \frac{\tui}{\ \newq\ }\ 
        \left [\ u_{xx} - \newq\ \right ]^2 dx}_{\I}\\ \\
\left. \begin{array}{l}\!\!\!
-\sumint\!\! \uxx
        \left [ \ \sWpp\ -\ \qq\ \right ]\!\!
 \\ \\
        \qquad\qquad\qquad\qquad\qquad\ \
        \times\ \left [\ \dif\ \right ]\ dx
\end{array}\right \} \mbox{\small{\II}}\\ \\
\underbrace{+\sm \utxi\ [\ \up\ -\ \ui\ ]}_{\CC}\\ \\
\underbrace{-\sumint \sWpp\ \uxx\ \q\ dx}_{\DD} ,
\end{array} \]
if for the values of $i$ such that \up\ =\ \um\ we define
\[ \qq \]
to be zero.  For such an $i$ the $i$th term in \I\ +\ \II\ is
\[ \quarto. \]
The term \I\ is analogous to the term\ \, $- \int_0^1 |u_t\ -\ u_t^m|^2\ dx$\/
in \cSection~3. It is negative.  The analogue of \II\ is identically
zero for the heat equation.  So \II\ is an error term due to the
nonlinearity of the equation.  The terms $C$\/ and $D$\/ will combine
much as $A\ +\ B\ +\ C$\/ did above.

        To estimate $C\ +\ D$, it is convenient to introduce \Wb,
\begin{equation} 
        \Wb^{\prime\prime}(y)\ \stackrel{\triangle}{=}
         \sqrt{1\ +\ y^2}\ \WD(y) .
\label{tildew}
\end{equation}
Notice that Eq.~($\thecurvature .1$) can be written in divergence form
\[
        u_t\ =\ (\tilde{W}^{\prime}(u_x))^{\prime} .
\]
We proceed essentially as we did for the heat equation.  On the one hand,
\begin{equation}\begin{array}{lcl}
        C&=&\half\ \s\ [\ \up\ -\ \ui\ ]\ \cdot \ 
                \left [\ \utxi\ \right ]\\
                &&\\ 
        &&\!\!\!\!\!\!\!+\ \half\ \s\ 
                [\ \ui\ -\ \um\ ]\ \cdot \ 
                \left [\ \utxm\ \right ]\ +\ \mbox{(v)}, 
\end{array}
\end{equation}
where of course
\begin{eqnarray*}
        \utxi\ &=&\ \iut\\ &&\\
        &=&\mbox{} \left. \frac{d}{dx}\ 
                \Wb^\prime(u_x(x,t))\ \right |_{x\ =\ x_i(t)} .
\end{eqnarray*}
On the other hand, $D$ equals
\[
-\s [\ \Wbi\ -\ \Wbm\ ]\ \q\ = \]
\[ \ \ \half \s\ [ \up\ -\ \ui ]
        \left [ \ -\ \frac{\Wbi\ -\ \Wbm\ }{\li\ } \right ]\ r_i\]\ \,
\[ +\ \half \s\ [ \ui\ -\ \um ]  
        \left [ \ -\ \frac{\Wbi\ -\ \Wbm\ }{\li\ } \right ]\ r_i\ .
\]
Therefore $C\ +\ D$\/ equals
\begin{equation}\begin{array}{l}
        \displaystyle   +\ \s\ \onefourth\ [\ \up\ -\ \ui\ ]\ \cdot \ 
                \left. \frac{d^2}{dx^2}\ \left [\ \mal\ \right ]
                \right |_{x\ =\ \xi_i(t)}\ r_i\ l_i
\\  \\  \displaystyle -\ \s\ \onefourth\ 
                [\ \ui\ -\ \um\ ]\ \cdot \ 
                \left. \frac{d^2}{dx^2}\ \left [\ \mal\ \right ]
                \right |_{x\ =\ \zeta_i(t)}\ r_i\ l_i
\\ \\   \displaystyle +\ \half\ u_t(x_1(t),t)\ 
                [\ \numxt\ -\ \umxo\ ]\ (1\ -\ r_1)\ +\
                \mbox{(vi)}.
\end{array}\label{PPP}\end{equation}
for some\ $\xi_i\ \mbox{and}\ \zeta_i\ $\ belonging to $(x_{i-1},\ 
x_i)$. But,
\[ 
        \half\ u_t(x_1(t),t)\ [\ \numxt\ -\ \umxo\ ]\ (1\ -\ r_1)\ =
\qquad\qquad\qquad\qquad\qquad\qquad\ \ \ \
\] 
\vspace{1 mm}
\[
        \ \half\ [\ \numxt\ -\ \umxo\ ]\ 
                \cdot \ 
                \left. \frac{d^2}{dx^2}\ \left [\ \mal\ \right ]
                \right |_{x\ =\ \varphi_1(t)}
                (1\ -\ r_1)\ l_1\ +\ \mbox{(vii)},
\]
for some $\varphi_1$ belonging to $(0, x_1)$.  For zero boundary
conditions each of (i) through (vii) vanishes.  The terms
corresponding to $\xi_1$\/ and $\varphi_1$\/ combine and we obtain
\begin{equation}
        |\ C\ +\ D\ |\ \leq\ \half\ m\ \left \|\ 
        \frac{d^2}{dx^2}\ \Wb^\prime(u_x)\ \right \|_{\infty}(t) .
\label{sumcandd}
\end{equation}

        It remains to estimate \II, which we write \II\ =\
$\sum_i$\ \IIi.  If $\up\ =\ \um$\/ then
\[
        \IIi\ =\ \quarto
\]
is negative.  So we need only consider $i$\/ such that $\up\ \neq\ \um$.
By the inequality $ab\ \leq \frac{\delta}{2} a^2\ +\ \frac{1}{2\delta} b^2$,
\[
        |\ \IIi\ |\ \leq\ 
\underbrace{
        \frac{\delta}{2}\ \integ \left [\ \dif\ \right ]^2\ dx
}_{\IIoi}
\] \[
\underbrace{
        +\ \frac{1}{2\ \!\! \delta}\ \integ \uxxs\ 
        \left [\ \sWpp\ -\ \qq\ \right ]^2\ dx ,
}_{\IIti}
\]
with $\delta$\/ an arbitrary positive number.  We shall control
$\IIoi$\/ later, by showing that it is dominated by the corresponding
term of \I, if $\delta$\/ is sufficiently small.
In order to estimate the term $\IIti$, we rewrite
the two terms on the right hand side of
Eqs.~(\ref{delta}) using Taylor's expansion,
\[ \begin{array}{cr}
        \!\!\!\!\pWi\!\!\!\!\!&=\ \!\dWi\ +\ \!\half\ \WDi\ [\ \difi\ ]
\\ &\\
        &       +\ \onesixth\ \Wti\ [\ \difi\ ]^2,\ \
\\ &\\ 
        \!\!\!\!\pWm\!\!\!\!\!&=\ \!\dWi\ -\ \!\half\ \WDi\ [\ \difm\ ]
\\ &\\
        &       +\ \onesixth\ \imWt\ [\ \difm\ ]^2,
\end{array}\]
where \tei,\ \ite\ belong to the interval with endpoints\ \ui\ 
and\ \up.  We have not yet used the discretized differential 
equation;  recall that $\tui\ =\ \suxis\ \frac{\Delta_i}{l_i}\ r_i$.
We get
\begin{eqnarray}
        \qq\ &=&\suxis\ \frac{\deltai}{\ \ \difct\ \ }\nonumber\\ 
                \nonumber&&\\
        &=& \suxis\ \WDi\ +\ \Ehat, \label{P}
\end{eqnarray}
where
\begin{equation}\begin{array}{lcl}
        |\ \Ehat\ |&\leq&\ \onethird\ \suxis\ |\ \Wti\ |\ |\ \difi\ |
\\ && \\ &&\!\!\!\!\!\!\!\ +\ \onethird\ \suxis\ |\ \imWt\ |\ |\ \difm\ |.
\end{array}\label{PP}\end{equation}
We use here the fact that $|\ \up\ -\ \ui\ |,\ |\ \ui\ -\ \um\ |\ \leq\ 
|\ \up\ -\ \um\ |$,\ since \ui\ is between \um\ and \up, because
$\up\ \neq\ \um$.  It follows that
\[
        \sWpp\ -\ \qq\ =\ \qquad\qquad\qquad\qquad
\] \[ 
\qquad\qquad\qquad\qquad\qquad\qquad\qquad\qquad\ \ \ 
\Wbpppz\ [\ u_x\ -\ \ui\ ]\ -\ \Ehat , \]
for some $\eta_i$ in the interval with endpoints
\ui\ and $u_x(x,t)$, and thus
\[\begin{array}{lcl}
        |\ \IIti\ |&\leq&\displaystyle \frac{1}{\delta}\ \integ \uxxs\
                [\Wbpppz]^2\ |\ u_x\ -\ \ui\ |^2\ dx \\
                &&\\
&&\displaystyle\underbrace{
        +\ \frac{1}{\delta}\ \integ \uxxs\ \Ehats\ dx
}_{\bE} .
\end{array}\]
These terms are clearly controllable.

        We reorganize, 
$\ I\ +\ \II\ +\ C\ +\ D\ $ = $\ (\ I\ +\ \IIo\ )\ +\ (\ \IIt\ -\ E\ )
\ +\ C\ +\ D\ +\ E$ (there is no risk of confusing this $E$
with the energy).  We have shown that
\[
        \ddt\ \half\ \inte |\ u_x\ -\ u_x^m\ |^2 \ dx\ =\ 
        \underbrace{\s\ \IIIi}_{\III}\ 
        +\ \underbrace{\s\ \IVi}_{\IV}\ 
        +\ C\ +\ D\ +\ \underbrace{\s E_i}_{\EEE}\ ,
\]
with
\[ \begin{array}{lclc} &&&\\
         \IIIi\ &\stackrel{\triangle}{=}&
                \left \{ 
                \begin{array}{l}
                 -\ \integ \xiiii\  \left [\ \dif\ \right ]^2\ dx 
\\ \qquad\qquad\qquad\qquad\qquad\qquad\qquad\qquad\ \ \,
\mbox{if\ \ \up\ $\neq$\ \um}\\
                \quarto
                        \ \,\ \ \ \mbox{if\ \ \up\ =\ \um}
                \end{array}\right.
\!\!\!\!\!&,\\ &&&\\
        \IVi\ &\stackrel{\triangle}{=}&
                \left \{ 
                \begin{array}{l}
        \displaystyle\frac{1}{\delta}\ \integ \uxxs\
                [\Wbpppz]^2\ |\ u_x\ -\ \ui\ |^2\ dx\ 
\\ \qquad\qquad\qquad\qquad\qquad\qquad\ \ 
\qquad\qquad
\mbox{if\ \ \up\ $\neq$\ \um}
\\0\qquad\qquad\qquad\qquad\qquad\qquad\qquad\qquad\, 
\mbox{if\ \ \up\ =\ \um}
                \end{array}\right.
\!\!\!\!\!&,\\ &&&\\
        E_i\ &\stackrel{\triangle}{=}&
                \left \{ 
                \begin{array}{l}
        \displaystyle\frac{1}{\delta}\ \integ \uxxs\ \Ehats\ dx 
                \qquad\qquad\qquad\qquad\ \ \,
                \mbox{if\ \ \up\ $\neq$\ \um}
                \\ \\0 
                        \qquad\qquad\qquad\qquad\qquad\qquad\,
\qquad\qquad
                \mbox{if\ \ \up\ =\ \um}
                \end{array}\right.
\!\!\!\!\!&,\\ &&&
\end{array}
\]
and
\[
        \xiiii\ \stackrel{\triangle}{=}\ \qq\ -\ \frac{\delta}{2}.
\]
We can now prove the following

\vspace{2mm}

\noindent {\bf Theorem:}\ {\em
        Let $u$\/ be a $C^3$\/ solution of\/~{\rm Eq.~($\thecurvature .1$)}
with homogeneous Dirichlet boundary condition\/~{\rm Eq.~($\thecurvature .2$)}.
Let $u^m$\/ solve the discretized problem presented in\/~{\rm \cSection~2}.
Suppose that the initial data for $u^m$\/ satisfies
\begin{equation}
        ||\ (u_0^m)_x\ -\ (u_0)_x\ ||_\infty\ \leq\ m 
\label{tone}
\end{equation}
and
\begin{equation}
        |\ [(u_0^m)_x]_i\ -\ [(u_0^m)_x]_{i-1}\ |\ \leq\ m\ \qquad
                        \mbox{for all\ }i .
\label{ttwo}
\end{equation}
Then
\begin{equation}
        \Upsilon\ \stackrel{\triangle}{=}\ \half\ \inte
                |\ u_x\ -\ u_x^m\ |^2 \ dx
\label{tthree}
\end{equation}
satisfies a differential inequality of the form
\begin{equation}
        \frac{d\Upsilon}{dt}\ \leq\ \Lambda(t)\ \Upsilon\ +\ \Gamma(t)\ m
\label{tfour}
\end{equation}
when $m$\/ is sufficiently small.  The coefficients $\Lambda$\/ and
$\Gamma$\/ are independent of $m$, though they depend upon the
underlying solution $u$.  In particular, for any $T\ >\ 0$,
\[
        \sup_{0\leq t\leq T}\ ||\ u\ -\ u^m\ ||_{H^1([0,1])}\ 
                                \leq\ C_T\ m^{\frac{1}{2}} .
\]
}

\vspace{2mm}

        {\em Proof.}\ \ Before we start, note that it is
always possible to choose initial data satisfying 
conditions~(\ref{tone}) and (\ref{ttwo}).  This was discussed
in \cSection~3.

        Recall from \cSection~2 that the discrete evolution never
introduces new faces after the initial time.  This fact together
with hypothesis~(\ref{tone}) implies that $\max_i \{ \ui \}$\/
is bounded with the bound uniform in $m$.  Furthermore,
because of (\thejump) on page~\pageref{back} and hypothesis~(\ref{ttwo}),
\[
        \max_i\ |\ \ui\ -\ \um\ |\ \leq\ m
\]
for all time, not just for time zero.

        Almost all the work for the theorem has already been done.
The principal task that remains is to choose the
parameter $\delta$\/ in such a way that $\Xi_i\ \geq\ 0$\/ whenever
$\up\ \neq\ \um$.  From  Eq.~(\ref{P}) and inequality~(\ref{PP}) 
we see that there exist constants $c_1$\/ and $c_2$\/ such that
\[
        \qq\ \geq\ c_1\ -\ c_2\, m, 
\]
with both these quantities strictly positive.  Their values depend only on the
maximum and minimum slopes of the initial data, and on the 
properties of the ``surface energy'' function $W$.  We restrict $m$\/
to $m\ \leq\ c_1 / (2 c_2)$\/ and take $\delta\ =\ c_1$.  Then
$\Xi_i\ \geq\ 0$\/ and \III\ is nonpositive.

        We may now estimate $C$, $D$, $E$, and \IV.  By 
inequality~(\ref{sumcandd}),
\[
        |\ C\ +\ D\ |\ \leq\ \half\ m\ \left \|\ 
        \frac{d^2}{dx^2}\ \Wb^\prime(u_x)\ \right \|_{\infty}(t);
\]
by inequality~(\ref{PP}), 
\[
        |\ E\ |\ \leq\ {c_1}^{-1}\ (c_2\, m)^2\ \inte u_{xx}^2\ dx ;
\]
and from the expression for \IVi,
\[
        |\ \IV\ |\ \leq\ {c_1}^{-1}\ ||\ u_{xx}^2\ ||_{\infty}\
        \max_{\eta}\ [\tilde{W}^{\prime\prime\prime}(\eta)]^2 \cdot\
        \inte |\ u_x\ -\ u_x^m\ |^2 \ dx.
\]
The maximum in the last formula is not taken over all values of
$\eta$, but rather over the range of possible values of $u_x$\/
and $u_x^m$;  these values lie in a bounded interval which
depends only on the initial data for $u$.  Clearly,
\[
        |\ C\ +\ D\ +\ E\ |\ \leq\ \Gamma(t)\ m,
\]
with the function $\Gamma$\/ independent of $m$.  The estimate
for \IV\ can be expressed as 
\[
        |\ \IV\ |\ \leq\ \Lambda(t)\ 
                \half\ \inte |\ u_x\ -\ u_x^m\ |^2 \ dx.
\]
This yields inequality~(\ref{tthree}), and a standard application
of Gronwall's inequality gives inequality~(\ref{tfour}).\hfill $\Box$

\vspace{\baselineskip}

        Eq.~(\thecurvature) has the form
\begin{equation}
        u_t\ =\ \psi(u_x)\ W^{\prime\prime}(u_x)\ u_{xx},
\label{psi}
\end{equation}
with $\psi(\nu)\ =\ (1\ +\ \nu^2)^{\frac{1}{2}}$.  We have made
no particular use of the form of $\psi$.  A similar convergence
theorem holds for Eq.~(\ref{psi}) with any $\psi$, continuously
differentiable, positive, and bounded away from zero.

        One final remark.  For the case of the heat equation it was clear that
$m\ =\ \max_i |\ \ui\ -\ \um\ |$\/ was the parameter controlling
the accuracy of the approximation.  In the general case 
one might have expected something different, for example
that one should take a
coarser discretization where $W^{\prime\prime}$\/ is smaller.
Our estimates do not support such an idea.
In fact, the argument at the end of \cSection~3 suggests that if one
chooses the admissible slopes so that the distance between
any two adjacent ones is a fixed constant, then one might
be able to prove quadratic convergence in $m$.  
\setcounter{transfer}{\theequation}

\section{Convergence for the Neumann problem}

\setcounter{equation}{\thetransfer}
Here we study the Neumann problem
\begin{equation} \left \{ \begin{array}{l}\displaystyle
        \frac{u_t}{\suxs}\ =\ 
                \WD(u_x)\ u_{xx}\qquad \mbox{if} \ \ 0\leq x \leq 1\\
                u(x,0)\ =\ u_0(x) \\
        u_x(0,t)\ =\ a,\qquad
u_x(1,t)\ =\ b \end{array} \right. ,\label{neumann}\end{equation}
with $a$\/ and $b$\/ constants.
We consider only discretizations for which
$a$\/ and $b$\/ are admissible slopes and we choose
\[
        \umxo\ =\ a\qquad\mbox{and}\qquad \umxn\ =\ b
\]
for the boundary conditions to be satisfied.  

        We want to insure that faces $1$\/ and $N$\/ do not disappear
and that statement~(\thejump) holds.  So at each time we extend
$u^m$\/ in such a way that \umxoo\ and \umxo\ correspond to adjacent
corners of \Wt, \umxo\ lies between \umxoo\ and \numxt, and the length 
of the first face of the extension of $u^m$\/ is $l_1\ /\ r_1$\/ where
\addtocounter{equation}{+1}
\setcounter{ratio}{\theequation}
$$
        r_1\ =\ \frac{\ \numxt\ -\ \umxo\ }{\ \numxt\ -\ \umxoo\ }
\eqno(\theratio .1)
$$
(note that $0\ <\ r_1\ <\ 1$\/ and if the value of \numxt\ changes
then $r_1$\/ might jump).  
The reason for this choice of $r_1$\/ will become clear below,
when we prove convergence of the approximation scheme.
Similar conditions apply to the extension 
at the right endpoint with
$$
        r_N\ =\ \frac{\ \umxn\ -\ \umxnn\ }{\ \umxnp\ -\ \umxnn\ } .
\eqno(\theratio .2)
$$
Furthermore, we impose that for the 
extension of $u^m$\/ Eqs.~(\ref{delta}) and (\ref{K}) 
(for \deltai\ and \tui,\ respectively) are valid for $i$ equal to $1$\/
and $N$, i.e.\addtocounter{ldot}{-2}
$$
\begin{array}{lcl}
        \umto&=&\displaystyle\suxos\ \frac{\deltao}{\lo}\ r_1 ,\\
&&\\
        \umtn&=&\displaystyle\suxns\ \frac{\Delta_N}{l_N}\ r_N;
\end{array}
\eqno(\theldot ^{\prime})
$$
it is natural to do so since we want our scheme to amount to motion by 
weighted curvature.  The idea is that
the law governing the evolution of $u^m$ should be the same
in both the interior and boundary of the interval $[0, 1]$\/ 
(recall also, from the Remark on page~\pageref{rmk}, that for
homogeneous Dirichlet boundary conditions one could think of
$u^m$\/ as being defined on the whole real line, the extension
being odd with period two; in that case Eqs.~(\ref{delta}) and (\ref{K})
were satisfied for $i\ =\ 1$\/ and $i\ =\ N$, i.e.\ the law 
governing the evolution of $u^m$\/ was the same 
in both the interior and boundary of the interval $[0, 1]$).
\addtocounter{ldot}{+2}

        The function $u^m$ is determined by
solving Eqs.~$(\theldot .1^\prime)$\ and $(\theldot .2)$\ 
with the $c_i^j$'s given by Eqs.~(\ref{KK}), the $(u_t^m)_i$'s by 
Eqs.~\addtocounter{ldot}{-2}$(\theldot ^{\prime})$\addtocounter{ldot}{+2}
and (\ref{K}), the $\Delta_i$'s by Eqs.~(\ref{delta}), and $r_1$\/ and
$r_N$\/ by Eqs.~(\theratio); 
the initial data should satisfy condition~(\thekeep).

        Arguing as in \cSection~2, we see 
that faces $1$\/ and $N$\/ of the extension of $u^m$\/ and of $u^m$\/
do not disappear and that statement~(\thejump) holds, namely,
even as time evolves and certain faces disappear, jumps in
$\ui$\/ correspond to adjacent corners of $\bar{W}$.

        Suppose now that we want $u^m$\/ to approximate
the motion of $u$.  Then we assume, in addition, that 
the number of faces of $u_0^m$\/ is 
\begin{equation}
        N\ \leq\ \frac{c}{m} ,
\label{inverse}
\end{equation}
with $c$\/ a constant.  One can, for example, construct
$u_0^m$\/ using the method described on page~\pageref{constroi} 
(basically, by taking the union of segments on lines
tangent to $u_0$\/ and with admissible slopes)
to get $N\ \leq\ [\mbox{Total Variation of\ }u_0^{\prime}]\ / m$.

        We want to estimate the $H^1$\/ norm of $u\ -\ u^m$.  The values
of $r_1$\/ and $r_N$\/ have been chosen so that the sum of the boundary
terms (i)~through~(vii) of the previous \csection\ vanish (see the Remark
on page~\pageref{boundary}).  So the estimate of the $L^2$\/ norm of
$u_x\ -\ u_x^m$\/ given in the previous \csection\ remains valid.  
However, to control the $H^1$\/
norm of $u\ -\ u^m$\/ we need some additional information,
since we cannot use Poincar\'{e}'s inequality.
For any $g:\ [0,1]\ \mapsto\ \mbox{\bf R}$\/ we have
\begin{equation}
        \inte |\ g\ -\ \langle g\rangle\ |^2\ dx\ \leq\ \frac{1}{\pi^2}\
                        \inte g_x^2\ dx ,
\label{series}
\end{equation}
with $\langle g\rangle\ =\ \int_0^1 g\ dx$, so it suffices to control
$\langle u\ -\ u^m\rangle$.
The rate of change of the average of $u$\/ is
\[
\begin{array}{lcl}
        \ddt\ \inte u\ dx&=&\inte \ut\ dx\\
&&\\
                &=&\displaystyle \inte \frac{d}{dx}\ \Wtp(u_x)\ dx\\
&&\\
                &=&\Wtp(\uxn)\ -\ \Wtp(\uxo)\ =\ \Wtp(b)\ -\ \Wtp(a),
\end{array}
\]
with \Wb\ as in Eq.~(\ref{tildew}).
On the other hand, the rate of change of the average of $u^m$\/ is 
\[
        \ddt\ \inte u^m\ dx\ =\ \ddt\ \sumint u^m\ dx\ =\ \sumint \tui\ dx
\]
since $u^m$\/ is continuous and $\dot{x}_0\ =\ \dot{x}_N\ =\ 0$.  Hence,
by Eqs.~(\ref{P}) and (\ref{PP}),
\[
\begin{array}{lcl}
        \ddt\ \inte u^m\ dx&=&\s \suxis\ \deltai\ r_i\\
&&\\
                &=&\s \suxis\ \WDi\ \difct\ r_i\\
&&\\
                &&+\ \s \mbox{O}(m^2) .
\end{array}
\]
By assumption~(\ref{inverse}),
\[
\begin{array}{lcl}
        \ddt \inte u^m\ dx&=&\s \WbDi\ \difct\ r_i\ +\ \mbox{O}(m)\\
&&\\
                &=&\displaystyle\s \frac{\ \Wtp[\up]\ -\ \Wtp[\um]\ }{2}\ r_i
                        \ +\ \mbox{O}(m)\\
&&\\
                &=&\displaystyle
        \frac{\Wtp[\umxnp]r_N\ +\ \Wtp[\umxn]\ +\ \Wtp[\umxnn](1-r_N)}{2}\\
&&\\
        &&\displaystyle -\ 
        \frac{\ \Wtp[\numxt](1-r_1)\ +\ \Wtp[\umxo]\ +\ \Wtp[\umxoo]r_1}{2}
\\&&\\&&
        +\mbox{O}(m)\\
&&\\
        &=&\Wtp[\umxn]\ -\ \Wtp[\umxo]\ +\ \mbox{O}(m)\\
&&\\
        &=&\Wtp(b)\ -\ \Wtp(a)\ +\ \mbox{O}(m).
\end{array}
\]
(Observe that in this calculation
we only need that $0\ <\ r_1, r_N\ \leq\ 1$;  if the distance
between any two adjacent admissible slopes is a fixed constant then,
with the choice of $r_1$\/ and $r_N$\/ made in Eqs.~(\theratio),
one can improve this result to $\frac{d}{dt}\ \int_0^1 u^m\ dx\ =\ 
\Wtp(b)\ -\ \Wtp(a)\ +\ \mbox{O}(m^2)$.)
It follows that
\begin{equation}
        \left |\ \ddt\ \inte (u\ -\ u^m)\ dx\ \right |\ \leq\ C\, m .
\label{average}
\end{equation}
The constant $C$\/ depends only on $\max_i\, |\ \ui\ |$\/ and
on the constant appearing in hypothesis~(\ref{inverse}).
If the initial data satisfies inequality~(\ref{tone}) 
then $\max_i\, |\ \ui\ |$\/ is uniformly bounded.  If, in 
addition, the initial data satisfies inequality~(\ref{ttwo}) then,
by combining Eqs.~(\ref{tfour}), (\ref{series}), and (\ref{average}),
we easily deduce that
\[
        \sup_{0\leq t\leq T}\ ||\ u\ -\ u^m\ ||_{H^1([0,1])}\ 
                                \leq\ C_T\ m^{\frac{1}{2}} .
\]
\setcounter{transfer}{\theequation}

\section{Setup of the general Dirichlet problem}

\setcounter{equation}{\thetransfer}
In this \csection\ we show how one can set up a crystalline algorithm
for the general Dirichlet problem
\begin{equation} \left \{ \begin{array}{l}\displaystyle
        \frac{u_t}{\suxs}\ =\ 
                \WD(u_x)\ u_{xx}\qquad \mbox{if} \ \ 0\leq x \leq 1\\
                u(x,0)\ =\ u_0(x) \\
        u(0,t)\ =\ a(t),\qquad
u(1,t)\ =\ b(t)  \end{array} \right. .\label{dirichlet}\end{equation}
There is a new feature in this case: the time dependent boundary
condition can lead to the {\em creation of new faces}\/ at $x\ =\ 0$\/
and $x\ =\ 1$.  This causes our convergence argument to break down,
and, in fact, we do not prove that the algorithm converges.
It seems natural to choose
\addtocounter{ldot}{-2}
$$
        \umto\ =\ \ap\ 
\qquad\mbox{and}\qquad \umtn\ =\ b^{\prime}, 
\eqno(\theldot ^{\prime\prime})$$ 
instead of Eqs.~(\theldot), so the boundary conditions are
satisfied.\addtocounter{ldot}{+2}

        We also want to insure statement~(\thejump), i.e.\
that even as time evolves and certain faces disappear, jumps
in $\ui$\/ correspond to adjacent corners of $\bar{W}$.  To do this
we impose the following condition:
{\em there should be an extension of $u^m$\/ for which\/
{\rm Eqs.~(\ref{delta})} and\/~{\rm (\ref{K})}}\/ (for \deltai\ and \tui,\ respectively)\/
{\em are valid for $i$ equal to 1 and $N$,
and such that \umxoo\ and \umxo,\ and \umxnp\ and \umxn\ 
correspond to adjacent corners of \Wt.}\ 
(For the extension of $u^m$,
faces 1 and $N$ are now longer in general.  In the following we
focus our attention on the left endpoint and of course corresponding
statements hold for the right one.)  Observe that:
\begin{itemize}
\item Eq.~(\ref{K}) for $i\ =\ 1$ forces $\Delta_1$ and $\ap\ =\ \umto$\ 
to have the same sign. 
\item Since $l_1$ is smaller 
than the length of the first face of the extension of $u^m$, we have
\begin{equation}
        \lo\ \leq\ \suxos\ \ \frac{\deltao}{\umto}
\label{good}\end{equation}
if\ $\umto\ \neq\ 0$
(otherwise there is no restriction on \lo).\ 
\end{itemize}
These conditions might break down when \ap\ changes sign or when
equality holds in inequality~(\ref{good}).  More specifically:
\begin{itemize}
        \item {\em There might not exist a \umxoo\ adjacent
to \umxo\ such that $\Delta_1$ and \ap\ have the same sign.}\
Such a \umxoo\ does not exist if, and only if,
$\ap\ <\ 0$ and $\numxt\ >\ \umxo,$\ or
$\ap\ >\ 0$ and $\numxt\ <\ \umxo$.  (When it does exist it is unique.)
        \item {\em The length $l_1$ might not satisfy 
inequality\/~{\rm (\ref{good})}}.
\end{itemize}
Under these circumstances, we have to allow a new
face to appear at the boundary.  The introduction of the new face 
makes it possible to meet the requirements above;  it is necessary
in two cases:
\begin{description}
        \item [Case $(i)$] {\em If \ap\ is zero and about to become negative 
and $\numxt\ >\ \umxo,$\ or
\ap\ is zero and about to become positive and $\numxt\ <\ \umxo,$\ 
we introduce a new face with slope \numxt,\ now $(u_x^m)_3$}\ 
(see Figure~3). With the new ordering  $\umxo\ >\ \numxt$\ and 
$\umxo\ <\ \numxt$,\ respectively, enabling us to pick 
the new \umxoo\ uniquely.
        \item [Case $(ii)$] {\em If equality holds in 
inequality\/~{\rm (\ref{good})} and $\dot{l}_1$ is positive we introduce a new 
face with slope \umxoo,\ now \umxo}\ (see Figure~4).
Note that when a new face appears it has zero length so thereafter 
Eq.~(\ref{good}) is satisfied during some nonzero time interval.
\end{description}

\includegraphics{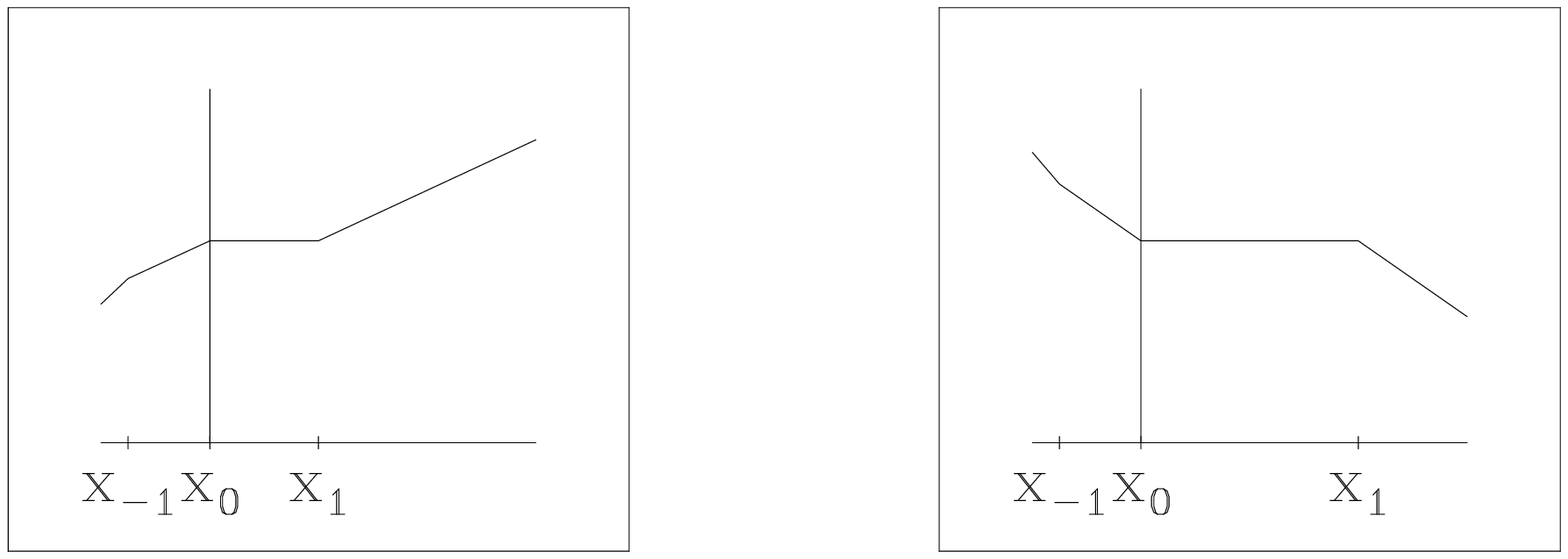}
\vglue 1.2truein                                              

\begin{center}
        Figure 3.\label{fthree} A new face is about to appear\\
        in the interval $[0, 1]$. Case $(i)$.
\end{center}

\newpage

\includegraphics{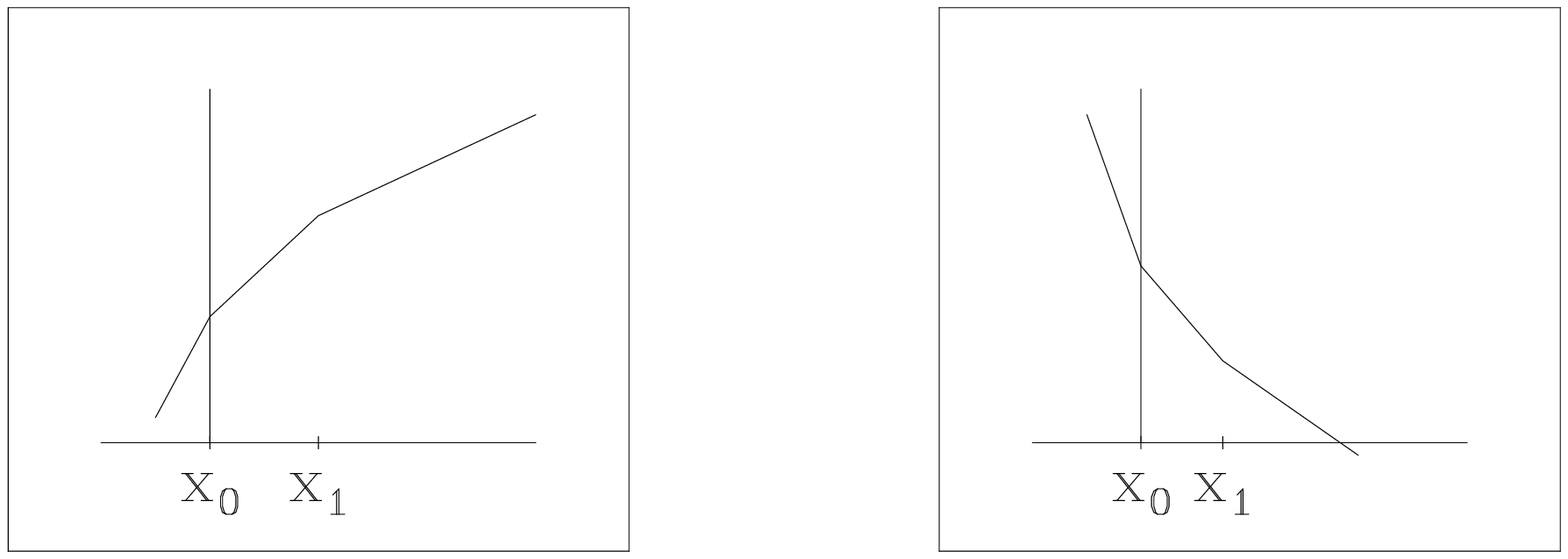}
\vglue 1.2truein

\begin{center}
        Figure 4.\label{ffour} A new face is about to appear\\
        in the interval $[0, 1]$. Case $(ii)$.
\end{center}

As long as neither Case $(i)$ nor Case $(ii)$ occurs
$u^m$ is determined by
solving Eqs.~$(\theldot .1^\prime)$\ and $(\theldot .2)$\ 
with the $c_i^j$'s given by Eqs.~(\ref{KK}), the $(u_t^m)_i$'s by 
Eqs.~\addtocounter{ldot}{-2}$(\theldot ^{\prime\prime})$\addtocounter{ldot}{+2}
and (\ref{K}), 
and the $\Delta_i$'s
by Eqs.~(\ref{delta}); the initial data should satisfy condition~(\thekeep).

        For general Dirichlet boundary conditions we have not been
able to prove convergence.  Nevertheless we examine how one 
can adapt part of the argument given in \cSection~4 to the present
case and see where it breaks down.
Here (contrary to what happened in \cSection~4) $(u_t^m)_1$\/ does not equal
$\sqrt{1\ +\ (u_x^m)_1^2\ }\ \frac{\Delta_1}{l_1}$, in general.  
For equality to hold we should substitute $l_1$\/ by the length of the 
first face of the extension of $u^m$ 
($\sqrt{1\ +\ (u_x^m)_1^2\ }\ \Delta_1\ /\ \umto$ if \umto\ is not zero,
otherwise it is not determined uniquely).  
We take 
\[
        r_1\ =\ \frac{\lo\ \umto}{\ \suxos\ \deltao\ }\qquad
\mbox{if}\ \ \umto\ \neq\ 0.
\]
We can take $r_1$ to be one if $\umto\ =\ 0$, since
$\umto\ =\ \sqrt{1\ +\ (u_x^m)_1^2\ }\ \frac{\Delta_1}{l_1}$, 
as $\Delta_1\ =\ 0$. Note that for the extension of $u^m$ 
we considered in \cSection~4
(odd and periodic) the ratio of the length of the first face of $u^m$\/
to the length of the first face of the extension of $u^m$ is one half,
but we could have extended $u^m$ so that this ratio was one.
Eq.~(\ref{good}) says that $r_1\ \leq\ 1$.  The length of the
first face of the extension of $u^m$\/ is $l_1\ /\ r_1$.

        With the $r_i$'s as in the previous paragraph we can estimate 
the growth of the $H^1$ norm of $u\ -\ u^m$
as was done in \cSection~4.  Now the boundary terms do not 
vanish,\addtocounter{ldot}{-2} but using 
Eqs.~$(\theldot ^{\prime\prime})$\addtocounter{ldot}{+2} we see 
that they add up to zero (see the Remark on page~\pageref{boundary}).
However, our argument also used  
that $\max_{1\leq i\leq N}\ |\ \ui\ |$\/ stays bounded as 
$m \rightarrow 0$, for $0\ \leq\ t\ \leq\ T$.  It seems plausible that
this should be the case if the solution of the differential 
equation~(\ref{dirichlet}) is such that 
$\sup_{0\leq t\leq T}\, ||\ u_x\ ||_{\infty}$\/ is finite.  This 
quantity is finite if, for example, we assume that
\addtocounter{equation}{+1}
\setcounter{growth}{\theequation}
$$
\frac{\cone}{\ouxs}\ \leq\ 
\suxs\ W^{\prime\prime}(u_x)\ \leq\ \frac{\ctwo}{\ouxs}
\eqno(\thegrowth .1)
$$
and 
$$
\suxs\ \ |\ W^{\prime\prime\prime}(u_x)\ |\ \leq\ \frac{\cthree}{\suxsc}
\eqno(\thegrowth .2)
$$
for some \cone, \ctwo, and \cthree\ $>\ 0$.
Condition~$(\thegrowth .1)$\/ assures that Eq.~(\ref{dirichlet})
is uniformly parabolic.  For a proof that $(\thegrowth .2)$\/ implies
$\sup_{0\leq t \leq T}\ ||\ u_x\ ||_\infty$\/ finite see
inequalities~($\mbox{V}\!\mbox{I}.5.10$) and ($\mbox{V}\!\mbox{I}.5.11$) 
in \lsu\ \cite{L} with their parameter $m$ equal to zero.
[Their $m$ is of course unrelated ours.  By giving
it values other than zero, one can get alternatives to 
inequalities~(\thegrowth).\
For example, in the case of the heat
equation the inequalities in \cite{L}\ 
just mentioned are satisfied with their $m$ equal to two instead of zero,
since in this case $W^{\prime\prime}(u_x)\ =\ 1\ /\ \suxs$.]
Conditions~(\thegrowth)\ hold for motion by weighted curvature,
as we verify in the Appendix.

        The proof of \cSection~4 would go through if one could bound 
$\max_{1\leq i\leq N}$\ $|\ \ui\ |$\/ for $0\ \leq\ t\ \leq\ T$\/
by a constant independent independent of $m$.
\setcounter{transfer}{\theequation}

\vspace{\baselineskip}

\section{Appendix:  Physical and mathematical context}

\vspace{\baselineskip}

\setcounter{equation}{\thetransfer}
We summarize here the relation between this work and
the literature on surface energy driven motion of
phase boundaries, especially the papers by Angenent and Gurtin~\cite{AG}
and Taylor~\cite{TW}.
Consider an interface between two phases moving isothermally according to the
balance of capillary forces and constitutive equations
compatible with thermodynamics. An evolution equation
for the interface is derived in \cite{AG}.
When the phases have the same energy 
and the kinetic coefficient (which measures the drag opposing
interfacial motion) is one, it has the form
(see Eq.~(4.11) of \cite{AG})
\begin{equation}
        V\ =\ [\ \ft\ +\ \fppt\ ]\ K,
\label{eqv}
\end{equation}
where $V$ is the normal velocity of the smooth interface and $K$ is its
curvature, $\theta$ is the angle from a fixed coordinate axis
to the normal to the interface, and $f$\/ (assumed smooth) 
is the interfacial energy per unit length.
So fix a coordinate system.  Let $u(\ \cdot\ ,t)$ be the interface
at time $t$ and 
\[
        \theta\ =\ \argt ,
\]
$0\ \leq\ \theta\ \leq\ \pi$,
be the angle between the normal (to the graph of $u$ with positive
coordinate in $y$) and the $x$-axis.
(We do not use the standard definition of $\arctan$
but rather one with range in the interval $[-\pi, 0]$.)  The expressions
of $V$ and $K$ in terms of $u$ are
\[
        V\ =\ \frac{u_t}{\suxs}\qquad \mbox{and}\qquad 
        K\ =\ \frac{u_{xx}}{\suxsc}.
\]
We are led to consider Eq.~(\thecurvature) by taking
\begin{equation}
        W(u_x)\ =\  f \argtb \ \suxs ,
\label{W}
\end{equation}
because
\begin{equation}
        W^{\prime\prime}(u_x)\ =\ \frac{\ \ft\ +\ \fppt\ }{\suxsc}.
\label{Wpp}
\end{equation}
The function $f$ is the energy per unit length (of
the interface) whereas $W$ is energy per unit length of the projection
of the interface on the $x$-axis.
Note that the right hand side of Eq.~(\ref{eqv}) is the negative
of the gradient of 
\[
        E(u)\ =\ \int W(u_x)\ dx\ =\ 
                \int f\argtb\ \suxs\ dx .
\]
Taylor calls the negative of the gradient of $E$\/
the weighted curvature of the interface (see Sections 2.2
and 2.3 of \cite{TW}).

The interfacial energy per unit length, $f$, is usually represented
in a polar diagram for $1/f$, called Frank diagram.
The function $f$
is said to be strictly stable when $f\ +\ f^{\prime\prime}\ >\ 0$.  This
condition corresponds to a strictly convex Frank diagram,
and to $W^{\prime\prime}\ >\ 0$.
On the other hand, $f$ is usually said to be crystalline if its convexified
Frank diagram is a polygon, 
and if the vertices of this polygon form the complete
set of globally convex sections of the Frank diagram (i.e.\
the original diagram and the convexified one meet only at the
vertices of the latter) (see Section~10.3 of \cite{AG}).  
We prefer a slightly broader definition.  We shall call
such an energy {\em strictly}\/ crystalline
and do not require that a crystalline energy satisfy the second 
condition.  Hence, if the Frank diagram is a polygon then $f$\/
is crystalline, but not strictly crystalline.  In general,
we denote by $\bar{f}$\/ the function whose
Frank diagram is the convexification of the Frank diagram of $f$.

        There is an equivalent characterization of crystalline energies.
It is obtained as follows.  
Let $n(\theta)\ \stackrel{\triangle}{=}\ (\cos\theta , \sin\theta)$.
For any surface energy $f$, we may
extend $f$ to $\mbox{\bf R}^2$ as a homogeneous function of degree one,
\[
        f_0 (x)\ \stackrel{\triangle}{=}\ \left \{ \begin{array}{ll}
        ||x||\ f\left ( \arg\ \displaystyle \frac{x}{||x||}\right ) 
                                        &\ \ \mbox{if}\ \ x\ \neq\ 0\\
                             0          &\ \ \mbox{if}\ \ x\ =\ 0
                        \end{array} \right . .
\]
Recall that the Fenchel transform of $f_0$ is the function 
$f_0^\ast :\ \mbox{\bf R}^2\ \longmapsto\ [-\infty, +\infty]$ given by
\[
        f_0^\ast (y)\ \stackrel{\triangle}{=}\ \sup_{x\in\mbox{\bf R}^2}
                            \ \{x\cdot y\ -\ f_0(x)\},
\]
and that the Wulff set of $f$ is 
\[
        W_f\ \stackrel{\triangle}{=}\ \{\ x\in \mbox{\bf R}^2\ |\
                x\cdot n(\theta)\ \leq\ f(\theta)\ \mbox{for all}\
                \theta\ \}.
\]
The function $f$ is crystalline if $W_f$\/ is polyhedral,
and it is strictly crystalline if in addition 
$f(\theta)\ >\ f_0^{\ast\ast}(n(\theta))\ \ (=\ \sup_{y\in W_f}
\ \{y\cdot n(\theta)\})$\/ unless $n(\theta)$\/ is normal to $\partial W_f$\/
(see, for example, Fonseca \cite{F}).  One can check that
$\bar{f}(\theta)\ =\ f_0^{\ast\ast}(n(\theta))$.
 
        An evolution equation for an interface with a strictly
crystalline energy is derived in \cite{AG} using the
same physical laws which gave Eq.~(\ref{eqv}).  For such an
energy Eq.~(\ref{eqv}) is backward-parabolic.  Therefore, one restricts
crystalline interfaces to a space $\cal{M}$, consisting of
continuous piecewise linear functions such that the normal 
to each face makes an angle with the $x$-axis
corresponding to one of the corners of the polygon $1 / \bar{f}\ \times\ n$.
The normal velocity of a face is (see Eq.~(10.12) of \cite{AG})
\begin{equation} 
        V\ =\ \sigma\ \frac{\tilde{\Delta}}{L},
\label{eqvd}
\end{equation}
where $\tilde{\Delta}$ and $\sigma$ are constants for each face, and
$L$ is the length of the face.  For a face with normal 
$n(\theta)$
\[
        \tilde{\Delta}\ =\ 
        [\ \bar{f}^\prime(\theta +0)\ -\ \bar{f}^\prime(\theta -0)\ ],
\]
and $\sigma\ =\ +1$ if $\theta$
increases across the face, $\sigma\ =\ -1$ if $\theta$
decreases across the face, and $\sigma\ =\ 0$ if 
the face is nontransitional.
Geometrically, $\tilde{\Delta}$ is the length of the segment in
the Wulff set of $f$
with normal $n(\theta)$. Using this formula, one can 
compute the velocity 
$V_i$ in terms of $f$, the
$\theta_i$'s, and $L_i$
(the subindex $i$ refers to the $i$th face):
\begin{equation}
\begin{array}{lcl}
        V_i&=&\displaystyle
+\ \frac{1}{\bli}\ f(\teim)\ \csc (\tei\ -\ \teim)\\&&\\
           & &\displaystyle
-\ \frac{1}{\bli}\ f(\tei)\ [\ \cot (\tei\ -\ \teim)\ 
                                +\ \cot(\teip\ -\ \tei)\ ]\\&&\\
           & &\displaystyle
+\ \frac{1}{\bli}\ f(\teip)\ \csc (\teip\ -\ \tei).
\end{array}
\label{vifornext}
\end{equation}
Since
\[
        \ui\ =\ -\ \cot \tei,
\] \[
        V_i\ =\ \frac{\tui}{\sqrt{1 \ +\  \uti\ }},
\] \[
  W[\ui]\ =\ f \argtbi\ 
        \suxis\ =\ \frac{f(\tei)}{\sin\tei},
\] 
and $l_i\ =\ L_i \sin \tei$ (we recall that $l_i$ is the length of 
the projection of the $i$th face on the $x$-axis), one easily checks
that this is Eq.~(\ref{K}), i.e.\ 
$\frac{\Delta}{l}\ =\ \sigma\ \frac{\tilde{\Delta}}{L}$. Notice that
\[
        \bar{W}(u_x)\ =\ \bar{f}\argtb\ \suxs,
\] 
where $\bar{f}$\/ is as above and $\bar{W}$\/ is as in
\cSection~2.  Hence the right hand side of Eq.~(\ref{eqvd}) is the
negative of the gradient of
\[
        \bar{E}(u)\ =\ \int \bar{W}(u_x)\ dx\ =\ 
                \int \bar{f}\argtb\ \suxs\ dx
\]
on $\cal{M}$, in other words the weighted curvature of the interface 
(see Sections 4.2 and 4.3 of \cite{TW}).

        In summary, we have checked that if $W$ is given by Eq.~(\ref{W}),
Eq.~(\thecurvature) is
Eq.~(\ref{eqv}) (with initial data and boundary conditions), 
and Eq.~(\ref{K}) is Eq.~(\ref{eqvd}). 
Our approximation scheme is crystalline in the sense that it approximates
the motion of an interface with a strictly convex energy by
the motion of an interface with a strictly crystalline energy 
and in the sense that $\bar{f}$\/ is crystalline.  

        In \cSections~4 and 5, we proved convergence of the crystalline
approximation scheme for homogeneous Dirichlet and
Neumann boundary conditions, respectively.  In such cases it follows
from the maximum principle that $||\ u_x\ ||_{\infty}$\/ is bounded.
In \cSection~6 we discussed the general Dirichlet problem, and noted
that growth conditions of the form (\thegrowth) are sufficient
to prove $L^{\infty}$\/ bounds on $u_x$.  Let us verify
that these conditions hold when $W$\/ is determined by
Eq.~(\ref{Wpp}), with $f$\/ strictly stable and smooth.
Then condition~($\thegrowth .1$)\/ is obvious and 
condition~($\thegrowth .2$) follows from the identity
\[
        W^{\prime\prime\prime}(u_x)\ =\
        \frac{\ -\ 3 \xu \ft\ +\ \fpt\ -\ 3 \xu \fppt\ +\ \fpppt\ }{\suxsf}.
\]
\setcounter{transfer}{\theequation}

\end{document}